\documentclass{ACmod}

\DeclareFontFamily{U}{rsf}{}
\DeclareFontShape{U}{rsf}{m}{n}{
  <5> <6> rsfs5 <7> <8> <9> rsfs7 <10-> rsfs10}{}
\DeclareMathAlphabet{\mathscr}{U}{rsf}{m}{n}

\DeclareMathAlphabet{\mathgth}{U}{euf}{m}{n}

\DeclareFontFamily{U}{cyr}{}
\DeclareFontShape{U}{cyr}{m}{n}{
  <5> wncyr5 <6> wncyr6 <7> wncyr7 <8> wncyr8 <9> wncyr9 <10-> wncyr10}{}
\DeclareMathAlphabet{\mathcyr}{U}{cyr}{m}{n}

\input cyracc.def

\DeclareSymbolFont{bbold}{U}{bbold}{m}{n}
\DeclareSymbolFontAlphabet{\mathbbold}{bbold}

\makeatletter
\def\operator@font{\sf}
\makeatother

\newcommand{\cA}{{\mathscr A}}

\newcommand{\CP}{{\mathbb{CP}}}

\newcommand{\cC}{{\mathscr C}}
\newcommand{\cE}{{\mathscr E}}
\newcommand{\cF}{{\mathscr F}}

\newcommand{\cO}{{\mathscr O}}
\newcommand{\cP}{{\mathscr P}}

\newcommand{\cS}{{\mathscr S}}

\newcommand{\cU}{{\mathscr U}}

\newcommand{\cV}{{\mathscr V}}

\newcommand{\MF}{{\sf MF}}

\newcommand{\shu}{{\sf sh}}

\newcommand{\Eu}{{\sf Eu}}

\newcommand{\hres}{{\mathsf{hres}}}

\newcommand{\bone}{{\mathbf 1}}

\DeclareMathOperator{\Tr}{Tr}

\DeclareMathOperator{\Der}{Der}

\DeclareMathOperator{\wt}{wt}

\DeclareMathOperator{\End}{End}

\DeclareMathOperator{\SL}{SL}

\DeclareMathOperator{\imag}{Im}
\DeclareMathOperator{\id}{id}

\DeclareMathOperator{\KS}{KS}

\newcommand{\ra}{\rightarrow}

\newcommand{\C}{\mathbb{C}}

\newcommand{\Q}{\mathbb{Q}}
\newcommand{\Z}{\mathbb{Z}}

\renewcommand{\phi}{\varphi}

\usepackage{tikz}
\usepackage{amscd}
\usepackage[tocflat]{tocstyle}
\usetikzlibrary{positioning}
\numberwithin{equation}{section}

\begin{document}

\title{Categorical Saito theory, II: Landau-Ginzburg orbifolds}

\author[Junwu Tu]{%
Junwu Tu}

\address{Junwu Tu, Institute of Mathematical Science, ShanghaiTech University, Shanghai 201210, China.}
\email{tujw@shanghaitech.edu.cn}

\begin{abstract}
{\sc Abstract:} Let $W\in \C[x_1,\cdots,x_N]$ be an invertible polynomial with an isolated singularity at origin, and let $G\subset {{\sf SL}}_N\cap (\C^*)^N$ be a finite diagonal and special linear symmetry group of $W$. In this paper, we use the category $\MF_G(W)$ of $G$-equivariant matrix factorizations and its associated VSHS to construct a $G$-equivariant version of Saito's theory of primitive forms. We prove there exists a canonical categorical primitive form of $\MF_G(W)$ characterized by $G_W^{{\sf max}}$-equivariance. Conjecturally, this $G$-equivariant Saito theory is equivalent to the genus zero part of the FJRW theory under LG/LG mirror symmetry. In the marginal deformation direction, we verify this for the FJRW theory of $\big(\frac{1}{5}(x_1^5+\cdots+x_5^5),\Z/5\Z\big)$ with its mirror dual B-model Landau-Ginzburg orbifold $\big(\frac{1}{5}(x_1^5+\cdots+x_5^5), (\Z/5\Z)^4\big)$. In the case of the Quintic family $\mathcal{W}=\frac{1}{5}(x_1^5+\cdots+x_5^5)-\psi x_1x_2x_3x_4x_5$, we also prove a comparison result of B-model VSHS's conjectured by Ganatra-Perutz-Sheridan~\cite{GPS}.
\end{abstract}

 \maketitle

\setcounter{tocdepth}{1}

\section{Introduction}

\paragraph{{\bf Backgrounds and Motivations.}} This is a sequel to our previous work~\cite{Tu} devoted to define and study what might be called a "$G$-equivariant Saito's theory of primitive forms".  Indeed, in {\em loc. cit.} we realized Saito's theory through Barannikov's notion~\cite{Bar} of Variation of Semi-infinite Hodge Structures (VSHS) associated with the category of matrix factorizations. Thus, to define a $G$-equivariant version of it, we simply consider the VSHS's associated with the category of $G$-equivariant matrix factorizations.

The motivation to develop such a theory is, to the author's point of view, quite substantial. First of all, for an invertible polynomial $W$, it is proved in~\cite{HLSW}~\cite{LLSS} that Saito's theory of $W$  is mirror dual to the FJRW theory~\cite{FJRW} of $(W^t, G^{t,{\sf max}})$, i.e. the dual polynomial with its maximal diagonal symmetry group. But the FJRW theory is defined with all the subgroups of $G^{t,{\sf max}}$ containing the diagonal symmetry group $J^t$. To make sense of mirror symmetry for $(W^t, G^t)$ with $J^t\subset G^t\subset G^{t,{\sf max}}$, one needs a $G$-equivariant Saito theory of $W$. A first case study was initiated by He-Li-Li~\cite{HLL}. Secondly, Saito's theory (when $G$ is trivial) is generically semisimple, while adding a non-trivial group $G$ may produce non-semisimple theory. For example, this includes Calabi-Yau hypersurfaces  through the so-called Landau-Ginzburg/Calabi-Yau correspondence.  If the genuz zero theory is not semi-simple, one can not define its higher genus theory using Givental-Teleman's construction. This motivates our third point of studying the category of $G$-equivariant matrix factorizations. Indeed, there is an alternative approach to extract categorical enumerative invariants developed in~\cite{Cos1}\cite{Cos2}\cite{CalTu} which is by construction all genera. A comparison result between that approach with the current VSHS approach is not yet known. We refer to~\cite[Conjecture 1.8]{CLT} for a precise formulation of this conjectural comparison.

\paragraph{{\bf A $G$-equivariant Saito theory.}} Let $W\in \C[x_1,\cdots,x_N]$ be an invertible polynomial with an isolated singularity at origin, and let $G\subset {{\sf SL}}_N\cap (\C^*)^N$ be a finite diagonal and special linear symmetry group of $W$. Consider the category $\MF_G(W)$ of $G$-equivariant matrix factorizations. Our construction of a $G$-equivariant Saito theory proceeds in two steps (done in Section~\ref{sec:hodge} and Section~\ref{sec:splitting} respectively):
\begin{itemize}
\item[1.] We first construct in Theorem~\ref{thm:construction-vshs} a primitive and polarized VSHS associated with $\MF_G(W)$. This VSHS provides us with a necessary framework to define Saito's notion of primitive forms. Note that in order to construct such a VSHS, a crucial ingredient is the Hodge-to-de-Rham degeneration property of the category $\MF_G(W)$. This was proved by Halpern-Leistner and Pomerleano in~\cite[Corollary 2.26]{HalPom}.
\item[2.] We then apply a bijection result (Theorem~\ref{thm:bijection2}) proved in~\cite{AT} to construct categorical primitive forms through a certain linear algebra data: splittings of the non-commutative Hodge filtration. For the category $\MF_G(W)$, we prove in Theorem~\ref{thm:canonical} there exists a canonical splitting of the non-commutative Hodge filtration. Through the aforementioned bijection, this yields a canonical choice of categorical primitive form of $\MF_G(W)$. 
\end{itemize} 
These two steps yields a VSHS $\cV^{\MF_G(W)}$ together with a canonical primitive form $\zeta^{\MF_G(W)}$. From this data, one obtains a Frobenius manifold by a standard construction~\cite{SaiTak}. Conjecturally, this Frobenius manifold is mirror dual to the genus zero part of the FJRW theory associated with the dual orbifold $(W^t, G^t)$. In~\cite{HLSW}), the conjecture is proved when $G=\{\id\}$ and $G^t=G^{t,{\sf max}}_{W^t}$.

We also remark that the step $1.$ above, i.e. the construction of VSHS's holds in much more general setup, as the degeneration result~\cite[Corollary 2.26]{HalPom} is proved for any $W$ with proper critical loci, and with an arbitrary finite symmetry group. Step $2.$ is the main reason we restrict ourselves to invertible polynomials.

\paragraph{{\bf Non-commutative calculations in the case of cubics.}} In Section~\ref{sec:cubic} we work out an example of $G$-equivariant Saito theory when $W=\frac{1}{3}(x^3_1+x_2^3+x_3^3)$ is the Fermat cubic and $G=\Z/3\Z$ the diagonal action. In this section, all computations are done non-commutatively in the sense that we work with Hochschild and cyclic chain complexes. The non-commutative calculation is needed to perform higher genus computations of categorical enumerative invariants as defined by Costello~\cite{Cos1}\cite{Cos2}, see also~\cite{CalTu}.

\paragraph{{\bf Calculations in the case of quintics.}} Section~\ref{sec:quintic} uses the comparison result in our previous work~\cite{Tu} to deduce  two main applications of the $G$-equivariant  Saito theory:
\begin{itemize}
\item In Theorem~\ref{thm:b-model-comparison}, we prove that the categorical VSHS is isomorphic to the geometric one constructed by Griffiths in the case of the mirror quintic family. The proof relies on Griffiths~\cite{Gri} and Carlson-Griffiths~\cite{CarGri}'s calculation of VSHS's related to projective hypersurfaces.  This $B$-model comparison result is needed in Ganatra-Perutz-Sheridan's work~\cite{GPS} which aims to establish a conceptual proof of the classical mirror symmetry conjecture via the categorical equivalence proved by Sheridan~\cite{She2}. In~\cite[Conjecture 1.14]{GPS} a general conjecture was stated about the comparison between categorical VSHS's and Griffiths' VSHS's. Theorem~\ref{thm:b-model-comparison} does not yield this conjecture in full generality, but suffices for  Ganatra-Perutz-Sheridan's application in the case of quintic mirror symmetry.

\item Theorem~\ref{thm:mirror} proves a LG/LG mirror symmetry result. Namely, for the Fermat quintic $W=\frac{1}{5}(x_1^5+\cdots+x_5^5)$, we match the equivariant Saito theory of $\MF_{(\Z/5\Z)^4}(W)$ (defined using the canonical primitive form $\zeta^{\MF_{(\Z/5\Z)^4}(W)}$) with the genus zero FJRW theory of $(W, \Z/5\Z)$, both in the marginal deformation direction. The proof is based on our previous comparison result~\cite{Tu} and a technique of computing primitive forms developped by Li-Li-Saito in~\cite{LLS}. Once the primitive form is calculated, the associated genus zero prepotential function is matched easily with Chiodo-Ruan's calculation of the FJRW invariants~\cite{ChiRua}.
\end{itemize}

\paragraph{{\bf Conventions.}} We following Sheridan's sign convention and notations in~\cite{She}. In particular, the notation $\mu_n (n\geq 1)$ stands for higher $A_\infty$ products in the shifted sign convention. For a sign $|a|\in \Z/2\Z$, denote by $|a|'=|a|+1\in \Z/2\Z$ for its shift. Furthermore, as in {{\sl loc. cit.}}, for a $\Z$-graded vector space/module $M$, the notation $M[[u]]$ stands for the $u$-adic completion (with $\deg(u)=-2$) of $M[u]$ in the category of $\Z$-graded vector spaces/modules. For example, the completion of the ground field $\mathbb{K}[[u]]=\mathbb{K}[u]$ has no effect.  

Following~\cite[Section 3]{She}, for a Hochschild cochain $\phi\in C^*(A)$ of an $A_\infty$-algebra, the notations $b^{1|1}(\phi)$, $B^{1|1}(\phi)$, $\iota(\phi)$, $L_\phi$ stands for certain actions of Hochschild cochains on Hochschild chains $C_*(A)$ or cyclic chains $C_*(A)[[u]]$, explicitly defined by 
\begin{align*}
b^{1|1}(\phi)(a_0|\cdots|a_n):= & \sum (-1)^{|\phi|'\cdot (|a_0|'+\cdots+|a_j|')} \mu\big( a_0,\ldots,\phi(a_{j+1},\ldots),\overbrace{\ldots,a_k}\big) | \ldots a_n,\\
B^{1|1}(\phi)(a_0|\cdots|a_n):= & \sum (-1)^{|\phi|'\cdot (|a_0|'+\cdots+|a_j|')} \bone|a_0|\ldots|\phi(a_{j+1},\ldots)|\overbrace{\ldots,a_n},\\
\iota(\phi):= & b^{1|1}(\phi)+uB^{1|1}(\phi),\\
 L_\phi(a_0|a_1|\cdots |a_n):=  & \sum (-1)^{|\phi|'\cdot (|a_0|'+\cdots+|a_j|')} a_0|a_1|\cdots|a_j|\phi(a_{j+1},\cdots,a_{j+l})|\cdots|a_n\\
 &+\sum \phi(\overbrace{a_0,\cdots,a_{l-1}})|a_l|\cdots |a_n
\end{align*}

\paragraph{{\bf Acknowledgment.}} The author is deeply indebted to Lino Amorim, Andrei C\u ald\u araru, Yunfeng Jiang, Si Li, Sasha Polishchuk, Yefeng Shen and Nick Sheridan for useful discussions. It was Si Li who suggested that requiring $G_W^{\sf max}$-equivariance should fix a canonical primitive form, which is precisely the content of Theorem~\ref{thm:canonical}. Sasha Polishchuk pointed out to me the reference~\cite{HalPom} on the Hodge-to-de-Rham degeneration property. Special thanks to Andrei C\u ald\u araru, Rahul Pandharipande, and Nick Sheridan as organizers of a beautiful workshop held at the Forschungsinstitut f\" ur Mathematik, ETH Z\" urich where the author had the great opportunity to present results of the paper. In particular, Theorem~\ref{thm:b-model-comparison} was inspired by Sheridan's talk at the workshop who also suggested a careful treatment of convergence issues in Section~\ref{sec:quintic}. 

\section{Non-commutative Hodge theory}~\label{sec:hodge}

In this section, we recall the categorical construction of VSHS's. The important notion of VSHS's was introduced by Barannikov~\cite{Bar}. For its categorical construction, our primary references are~\cite{She}\cite{Shk2}\cite{Shk3}\cite{Shk4}\cite{CLT}.

\paragraph{{\bf VSHS's.}} Let $(R,\mathfrak{m})$ be a complete regular local ring of finite type over $\mathbb{K}$. A polarized VSHS over $R$ is given by the following data:
\begin{itemize}
\item A $\Z/2\Z$-graded free $R[[u]]$-module $\mathcal{E}$ of finite rank.
\item A {{\sl flat}}, degree zero, meromorphic connection $\nabla$ such that it has a simple pole along $u=0$ in the $R$-direction, and at most an order two pole at $u=0$ in the $u$-direction. In other words, we have
\[ \nabla_X: \cE \ra u^{-1} \cE, \; X\in \Der(R), \; \;\; \nabla_{\frac{\partial}{\partial u}}: \cE \ra u^{-2} \cE.\]
\item A $R$-linear, $u$-sesquilinear, and $\nabla$-constant pairing $\langle-,-\rangle_{{\sf hres}}: \mathcal{E}\otimes \mathcal{E} \rightarrow R[[u]]$, such that its restriction to $u=0$ is non-degenerate.
\end{itemize}
We refer to~\cite[Section 2]{Tu} for more details.

\begin{Definition}
A VSHS $\big(\cE,\nabla,\langle-,-\rangle\big)$ over a regular local ring $R=\mathbb{K}[[t_1,\cdots,t_\mu]]$ is called primitive if there exists a section $\zeta\in \cE$ such that the map $\rho^\zeta: \Der_R \ra \cE/u\cE$ defined by
\[\rho^\zeta( \frac{\partial}{\partial t_j}) := \pi\big( u\nabla_{\frac{\partial}{\partial t_j}} \zeta \big)\]
is an isomorphism. Here $\pi: \cE\ra \cE/u\cE$ is the canonical projection map.
\end{Definition}

\medskip
\paragraph{{\bf Primitive VSHS's from Non-commutative geometry.}} Let $A$ be a $\Z/2\Z$-graded, finite-dimensional, smooth cyclic $A_\infty$ algebra of parity $N\in\Z/2\Z$. For $\epsilon\in \Z/2\Z$, denote by $HH_{[\epsilon]}(A)$ (respectively $HH^{[\epsilon]}(A)$) the parity $\epsilon$ part of the Hochschild homology (cohomology). 

The cyclic pairing on $A$ induces an isomorphism $A\ra A^\vee[-N]$ of $A_\infty$-bimodules, which yields an isomorphism
\[ HH^*(A) \cong HH^{*+N} (A^\vee).\]
The right hand side is naturally isomorphic to the shifted dual of the Hochschild homology $HH_{*-N}(A)^\vee$. There is also a degree zero pairing naturally defined on Hochschild homology, known as the Mukai-pairing
\[ \langle-,-\rangle_{{\sf Muk}}: HH_*(A)\otimes HH_*(A) \ra \mathbb{K}.\]
Since $A$ is smooth and proper, this pairing is non-degenerate. Thus, it induces an isomorphism $HH_{*}(A) \cong HH_{-*}(A)^\vee$. Putting together, we obtain a chain of isomorphisms:
\begin{equation}~\label{def:D}
D: HH^*(A) \cong HH^{*+N} (A^\vee)\cong HH_{*-N}(A)^\vee \cong HH_{N-*}(A)
\end{equation}
Let us denote by $\omega:=D(\bone)$, the image of the unit element $\bone\in HH^*(A)$. It was shown in~\cite{AT} that the duality map $D$ is an isomorphism of $HH^*(A)$-modules, i.e. it intertwines the cup product with the cap product.

Assume that $A$ satisfies the Hodge-to-de-Rham degeneration property. It was shown in~\cite{Iwa}\cite{Tu2} that the formal deformation theory of such $A$ is unobstructed. Let $\{\phi_1,\cdots,\phi_\mu\}$ be a basis of the even Hochschild cohomology $HH^{[0]}(A)$, and set
\[ R:=\mathbb{K}[[t_1,\ldots,t_\mu]].\]
Let $\cA$ denote the universal family of deformations of $A$ as a $\Z/2\Z$-graded unital $A_\infty$-algebra, parametrized by the ring $R$.  By construction, the Kodaira-Spencer map 
\begin{align*}
\KS: \Der_\mathbb{K}(R)  & \ra  HH^{[0]}(\cA),\\
\frac{\partial}{\partial t_j} & \mapsto [\frac{\partial}{\partial t_j}\mu^\cA]
\end{align*}
is an isomorphism. Here $\mu^\cA= \prod_n \mu_n^\cA$ is the $A_\infty$-structure maps of the family $\cA$.

In this section, we prove the following

\begin{Theorem}~\label{thm:construction-vshs}
Let the notations be as in the previous paragraph. Then the parity $N$ part of the negative cyclic homology $HC^-_{[N]}(\cA)$ of the universal family carries a  natural primitive VSHS over $R$.
\end{Theorem}

\begin{Proof}
Recall that the VSHS on $HC^-_{[N]}(\cA)$ was constructed by putting
\begin{itemize}
\item[--] in the $t$-direction the Getzler's connection~\cite{Get} 
\item[--] in the $u$-direction the connection operator defined in~\cite{KonSoi}\cite{KKP}\cite{Shk}\cite{CLT}
\item[--] the categorical higher residue pairing $\langle-,-\rangle_{{\sf hres}}$ defined in~\cite{Shk2}\cite{She}
\end{itemize}
The compatibility of the above data holds in general, see for example~\cite{CLT} and~\cite{She}. The non-degeneracy of the pairing is due to Shklyarov~\cite{Shk2}. What remains to be proved is the locally freeness of $HC^-_{[N]}(\cA)$ as a $R[[u]]$-module, and its primitivity.  We will first argue that the Hochschild cohomology $HH^{[0]}(\cA)$ is a locally free $R$-module. In~\cite{Tu2}, it was proved that the DGLA $C^*(A)[1]$ is homotopy abelian, i.e. there is an $L_\infty$ quasi-isomorphism
\[ \cU: HH^*(A)[1] \ra C^*(A)[1]\]
where the left hand side is endowed with the trivial DGLA structure. Tensoring with $R$ and pushing forward the universal Maurer-Cartan element $\beta:=\sum_{j=1}^\mu t_j \phi_j$ yields an $L_\infty$ quasi-isomorphism
\[ \cU^\beta: HH^*(A)[1]\otimes R \ra C^*(\cA)[1].\]
This proves that $HH^*(\cA)\cong HH^*(A)\otimes R$ which is clearly locally free $R$-module of finite rank. In particular, its even degree part $HH^{[0]}(\cA)$ is also locally free $R$-module of finite rank.

Next, we shall relate Hochschild homology with Hochschild cohomology via the duality isomorphism~\ref{def:D}. Indeed, it was proved in~\cite{Tu2} that the universal family $\cA$ can be chosen, after an gauge transformation, to be a cyclic universal family, with respect to the cyclic inner product $\langle-,-\rangle: A^{\otimes 2}\ra \mathbb{K}$ on $A$, extended $R$-linearly to $\cA$. Thus, there exists a duality isomorphism
\[ D: HH^{[0]}(\cA)\cong \big( HH_{[N]}(\cA) \big)^\vee \cong HH_{[N]}(\cA).\]
By the Nakayama lemma, $D$ is an isomorphism.

To this end, it remains to prove that $HC^-_{[N]}(\cA)$ is a locally free $R[[u]]$-module of finite rank. We may use the homological perturbation lemma for this. We may choose a homotopy retraction between
\[ HH_*(A)[[u]]\cong HC^-_*(A) \cong \big( C_*(A)[[u]], b+uB\big).\]
The first isomorphism exists as we assume the Hodge-to-de-Rham degeneration property. Tensoring with $R$ and $\mathfrak{m}$-adically complete the result yields a $R[[u]]$-linear homotopy retraction
\[ HH_*(A) \otimes R[[u]] \cong\big( C_*(A)[[u]], b+uB\big) \widehat{\otimes} R.\]
Here $\widehat{\otimes}$ stands for the $\mathfrak{m}$-adic completed tensor product. Adding the perturbation $L_{\cU_*\beta}$ to the right hand side yields exactly the complex $C_*(\cA)[[u]]$ which calculate the negative cyclic homology of $\cA$. After homological perturbation we obtain
\[ \big(  HH_*(A) \otimes R[[u]], \delta \big) \cong \big( C_*(A)[[u]]\widehat{\otimes} R, b+L_{\cU_*\beta}+uB\big).\]
We claim that the differential $\delta$ is equal to zero. Indeed, let us consider the $u$-filtration on both sides, and consider the following commutative diagram.
\[\begin{CD}
HH_*(A)\otimes R @>>>  HH_*(A) \otimes R[[u]]/u^{k+1} @>>>  HH_*(A) \otimes R[[u]]/u^k\\
@VVV  @VVV @VVV \\
HH_*(\cA) @>>> H^*\big(C_*(A)[[u]]\widehat{\otimes} R/u^{k+1}\big) @>>> H^*\big(C_*(A)[[u]]\widehat{\otimes} R/u^k\big)
\end{CD}\]
The point is that the homology of $C_*(\cA)$, i.e. $HH_*(\cA)$ is a locally free $R$-module of rank $\dim_\mathbb{K} HH_*(A)$ as we proved earlier, which implies that the induced differential $\delta$ on the upper left corner $HH_*(A)\otimes R$ is trivial. Induction on $k$ yields that the differential $\delta\equiv 0$.

The primitivity of this VSHS is obtained by choosing a primitive element
\[ \zeta:= \mbox{ any lift of } \; D(\bone),\]
with $D$ the duality map. Note that such a lift exists since $HC^-_{[N]}(\cA)$ is locally free.  To verify the primitivity of $\zeta$, we compute
\begin{align*}
\rho^\zeta(\frac{\partial}{\partial t_j}) &= \pi \big( u\nabla_{\frac{\partial}{\partial t_j}}\zeta\big)\\
&= -\KS(\frac{\partial}{\partial t_j})\cap D(\bone)\\
&= D(\KS(-\frac{\partial}{\partial t_j}))
\end{align*}
Since both $D$ and the Kodaira-Spencer map are isomorphisms, we deduce that $\rho^\zeta$ is also an isomorphism.
\end{Proof}

\begin{Lemma}
Let $HC^-_{[N]}(\cA)$ be the primitive VSHS over $R$ defined in Theorem~\ref{thm:construction-vshs}. Then there exists a unique vector field $\Eu\in {{\sf Der}}_\mathbb{K}(R)$ such that the differential operator $\nabla_{u\frac{d}{du}}+\nabla_\Eu$ is regular, i.e. it has no poles. 
\end{Lemma}

\begin{Proof}
Recall from~\cite{CLT} the $u$-connection is of the form
\[ \nabla_{u\frac{d}{du}}= u\frac{d}{du}+\frac{\Gamma}{2}+\frac{\iota(\prod_{n} (2-n)\mu_n)}{2u}.\]
Thus, the $u^{-1}$-component of the operator $\nabla_{u\frac{d}{du}}$ is given by capping with the Hochschild cohomology class
\[ [\frac{\prod_{n} (2-n)\mu_n}{2}] \in HH^*(\cA).\]
As such, if we set $\Eu:=\KS^{-1} [\frac{\prod_{n} (2-n)\mu_n}{2}]$, the combined operator $\nabla_{u\frac{d}{du}}+\nabla_\Eu$ is regular. If $\Eu'$ is another such vector field, we deduce that $\KS(\Eu-\Eu')=0$, but $\KS$ is an isomorphism, so $\Eu-\Eu'=0$.
\end{Proof}

\medskip
\begin{Definition}~\label{defi:primitive}
 An element $\zeta \in HC^-_*(\cA)$ is called a primitive form of the polarized VSHS $HC^-_*(\cA)$ if it satisfies the following conditions:
	\begin{itemize}
		\item[P1.] (Primitivity) The map defined by
		\[ \rho^\zeta: \Der (R) \ra  HC^-_*(\cA)/uHC^-_*(\cA), \;\; \rho^\zeta(v):=[u\cdot \nabla^{{\sf Get}}_v\zeta]\]
		is an isomorphism.
		\item[P2.] (Orthogonality) For any tangent vectors $v_1, v_2\in \Der (R)$, we have 
		\[\langle u\nabla^{{\sf Get}}_{v_1} \zeta, u\nabla^{{\sf Get}}_{v_2} \zeta\rangle \in R.\]
		\item[P3.] (Holonomicity) For any tangent vectors $v_1,v_2,v_3\in \Der (R)$, we have
		\[ \langle u\nabla^{{\sf Get}}_{v_1} u\nabla^{{\sf Get}}_{v_2} \zeta, u\nabla^{{\sf Get}}_{v_3} \zeta\rangle \in R\oplus u\cdot R.\]
		\item[P4.] (Homogeneity) There exists a constant $r\in \mathbb{K}$ such that
		\[ (\nabla_{u\frac{\partial}{\partial u}}+\nabla_{{\sf Eu}}^{{\sf Get}})\zeta= r\zeta.\]
	\end{itemize}
\end{Definition}

\paragraph{{\bf Splittings of the Hodge filtration.}} The definition of primitive forms is admitted quite complicated. We use an idea due to Saito~\cite{Sai1}\cite{Sai2} to construct primitive forms via certain linear algebra data known as splittings of the Hodge filtration. We first describe this notion more precisely. The filtration defined by
\[ F^k HC_*^-(A):= u^k\cdot HC_*^-(A), \;\; k\geq 0\]
is known as the non-commutative Hodge filtration on the negative cyclic homology of $A$. Denote by
\[ \pi: HC_*^-(A) \ra HH_*(A)\]
the natural projection map.

\begin{Definition}~\label{defi-splitting}
 A splitting of the Hodge filtration of $A$ is a linear map
\[ s: HH_{*}(A) \rightarrow HC^-_{*}(A),\]
such that it satisfies
\begin{itemize}
\item[S1.] {{\sl Splitting property.}} $\pi\circ s=\id$;
\item[S2.] {{\sl Lagrangian property.}} For any two classes $\alpha,\beta\in HH_{*}(A)$, we have $$\langle \alpha,\beta \rangle = \langle s(\alpha), s(\beta)\rangle_\hres.$$
\end{itemize}
\begin{itemize}
\item[S3.] {{\sl Homogeneity.}} A splitting $s$ is called a good splitting if $\nabla_{u^2\frac{\partial}{\partial u}} \imag (s) \subset \imag (s) \bigoplus u \cdot \imag (s)$.
\item[S4.] {{\sl $\omega$-Compatibility.}} Recall that $\omega=D(\bone)$ under the duality map~\ref{def:D}.  A splitting $s: HH_*(A) \ra HC_*^-(A)$ of the Hodge filtration is called  {{\sl $\omega$-compatible}} if $\nabla_{u\frac{\partial}{\partial u}} s(\omega) - r\cdot s(\omega)\in u^{-1} \imag (s)$ for some constant $r\in \mathbb{K}$.
\end{itemize}
\end{Definition}

\medskip
The following result can be used to construct categorical primitive forms. See~\cite{AT} for a detailed proof. Its commutative version was proved in~\cite{LLS}. The proof in the non-commutative version is almost identical to its commutative version.

\begin{Theorem}~\label{thm:bijection2}
	Let $HC^-_*(\cA)$ be a polarized primitive VSHS constructed as in Theorem~\ref{thm:construction-vshs}. Let $\omega=D(\bone)\in HH_*(A)$. Then there exists a natural bijection between the following two sets
	\begin{align*}
	\cP &:= \left\{ \zeta\in HC^-_*(\cA) | \;\zeta \mbox{ is a primitive form such that } \zeta|_{t=0,u=0}=\omega.\right\}\\
	\cS &:= \left\{ s: HH_{[N]}(A)\ra HC_{[N]}^-(A) | \; s \mbox{ is an $\omega$-compatible good splitting of parity $[N]$}.\right\}
	\end{align*}
	Note that in the second set $\cS$, the definition of a splitting of parity $[N]$ is the same as that of Definition~\ref{defi-splitting} by simply replacing $HH_*$ and $HC^-_*$ with $HH_{[N]}$ and $HC_{[N]}^-$. 
\end{Theorem}

\medskip

\paragraph{{\bf Saturated Calabi-Yau $A_\infty$-categories.}} Let $\cC$ be a saturated (i.e. compact, smooth and compactly generated) Calabi-Yau $A_\infty$-category. Assume also that $\cC$ satisfies the Hodge-to-de-Rham degeneration property. We may choose a split generator $E\in \cC$, and apply the construction of VSHS described above to the $A_\infty$-algebra ${{\sf End}}_\cC(E)$. By the Morita invariance of all the structures involved in the construction (see~\cite{She}), the resulting VSHS is independent of the generator $E$, up to isomorphism of VSHS's. For this reason, we shall denote by $$\cV^{\cC}:= \Big( HC^-_*\big( {{\sf End}}_\cC(E)\big),\nabla,\langle-,-\rangle_{\sf hres}\Big)$$ this VSHS, suppressing its dependence on $E$.

\section{Canonical splittings of invertible LG orbifolds}~\label{sec:splitting}

Let $W\in \mathbb{C}[x_1,\cdots,x_N]$ be an invertible polynomial. We refer to~\cite[Section 2.1]{HLSW} and~\cite[Section 1.2]{Kra} for the definition of invertible polynomials. Denote by 
\[ G_W^{{\sf max}}:=\{ (\lambda_1,\cdots,\lambda_N)\in (\mathbb{C}^*)^N\mid W(\lambda_1 x_1,\cdots,\lambda_Nx_N)=W.\}\]
the group of maximal diagonal symmetries of $W$. Fix a subgroup $G\subset G_W^{{\sf max}}\cap \SL_N$. In this section, we consider the category $\MF_G(W)$ of $G$-equivariant matrix factorizations of $W$. Since $G_W^{{\sf max}}$ is a commutative group, if $\phi$ is a $G$-equivariant morphism, then $g\phi$ remains a $G$-equivariant morphism for any $g\in G_W^{{\sf max}}$. Thus the maximal symmetry group $G_W^{{\sf max}}$ still acts on $\MF_G(W)$. The main result of this section is the following

\begin{Theorem}~\label{thm:canonical}
There exists a unique $G^{{\sf max}}_W$-equivariant good splitting in the sense of Definition~\ref{defi-splitting} of the non-commutative Hodge filtration of $\MF_G(W)$. Furthermore it is also compatible with the Calabi-Yau structure determined by $dx_1\wedge\cdots\wedge dx_N$.
\end{Theorem}

\medskip
\medskip
\begin{remark}
For the Hodge-to-de-Rham degeneration property of $\MF_G(W)$, Kaledin' Theorem~\cite{Kal} does not apply as the category $\MF_G(W)$ is only $\Z/2\Z$-graded. Nevertheless, this degeneration  is proved in~\cite[Corollary 2.26]{HalPom}. It is also known that the category $\MF_G(W)$ is a saturated, Calabi-Yau dg-category, see~\cite[Theorem 2.5.2]{PolVai}. Thus, we may apply Theorem~\ref{thm:construction-vshs} to obtain a primitive VSHS $\cV^{\MF_G(W)}$. The above Theorem combined with Theorem~\ref{thm:bijection2}, implies that there exists a canonical primitive form 
\[ \zeta^{\MF_G(W)}\in \cV^{\MF_G(W)}.\]
\end{remark}

\medskip
The proof of Theorem~\ref{thm:canonical} occupies the rest of the section.  The existence is proved as in~\cite{LLS}. Our contribution is the uniqueness. When $G$ is trivial, we would like to remark that the canonical splitting appears already in the works~\cite{HLSW}~\cite{LLSS} in the context of LG/LG mirror symmetry.


\paragraph{{\bf Reduction to commutative theory.}} Let us first assume that the group $G=\id$. Recall from our previous work~\cite{Tu} the set of splittings of  the non-commutative Hodge filtration of $\MF(W)$ is naturally in bijection with the set of splittings in the sense of Saito~\cite{Sai1}~\cite{Sai2}. Namely, one replaces
\begin{align*}
HH_*(\MF(W)) & \mapsto {\sf Jac} (W) dx_1\cdots dx_N\\
HC_*^-(\MF(W)) & \mapsto H^*\big( \Omega_{\C^N}^*[[u]], dW+ud_{DR}\big)\\
\nabla_{u\frac{\partial}{\partial u}} &\mapsto \nabla_{u\frac{\partial}{\partial u}}:= u\frac{\partial}{\partial u} - \frac{N}{2} -\frac{1}{u}\cdot W\\
\langle-,-\rangle_{{\sf hres}} & \mapsto \langle-,-\rangle_{{\sf hres}} \mbox{\; (higher residue pairing of Saito)}\\
\omega & \mapsto dx_1\cdots dx_N
\end{align*}

\paragraph{{\bf Weights and homological degrees.}} The key ingredient in proving Theorem~\ref{thm:canonical} is the existence of a rational grading. Recall that by definition of an invertible polynomial (see for example~\cite{HLSW}) that there exists a rational weights $\wt (x_j)=q_j, \; j=1,\ldots N$ such that $q_j\in \Q\cap (0,\frac{1}{2}]$. With respect to these weights, the polynomial $W$ is quasi-homogeneous of weight $1$, or equivalently we have
\begin{equation}~\label{eq:W}
W= \sum_{j=1}^N q_j x_j \frac{\partial W}{\partial x_j}
\end{equation}
We define a rational grading on the twisted de rham complex $ \big( \Omega_{\C^N}^*[[u]], dW+ud_{DR}\big) $ by
\[ \deg (x_i)  = -2 q_i, \;\; \deg (dx_i)  = -2 q_i +1, \;\; \deg (u) = -2.\]
With respect to this degree, the operator $dW+ud_{DR}$ is of degree $-1$.

\begin{Lemma}~\label{lem:omega}
The cohomology class $[dx_1\cdots dx_N]\in H^*\big( \Omega_{\C^N}^*[[u]], dW+ud_{DR}\big)$ satisfies 
\[ \nabla_{u\frac{\partial}{\partial u}} [dx_1\cdots dx_N]= -\frac{c_W}{2} [dx_1\cdots dx_N],\]
where the constant is $c_W:=\sum_{j=1}^N(1-2q_j)$ known as the central charge of $W$. 
\end{Lemma}

\begin{proof}
We compute using Equation~\ref{eq:W}:
\begin{align*}
 \nabla_{u\frac{\partial}{\partial u}} [dx_1\cdots dx_N] & =- \frac{N}{2} [dx_1\cdots dx_N] - \frac{1}{u} \sum_{j=1}^N q_j x_j \frac{\partial W}{\partial x_j} [dx_1\cdots dx_N]\\
 &=- \frac{N}{2} [dx_1\cdots dx_N] - \frac{1}{u} \sum_{j=1}^N q_j [x_j (-1)^{j-1}dW\wedge dx_1\cdots \widehat{x_j}\cdots dx_N]\\
 &= - \frac{N}{2} [dx_1\cdots dx_N] - \frac{1}{u} \sum_{j=1}^N (-1)^j q_j u \cdot [d_{DR} (x_jdx_1\cdots \widehat{x_j}\cdots dx_N)]\\
 & = -\frac{N}{2} [dx_1\cdots dx_N] + \frac{1}{u} \sum_{j=1}^N q_j u\cdot  [dx_1\cdots dx_N]\\
 & = -\frac{c_W}{2}  [dx_1\cdots dx_N].
 \end{align*}
\end{proof}

Let $f\in \mathbb{C}[x_1,\ldots,x_N]$ be a quasi-homogeneous polynomial. Computing $ \nabla_{u\frac{\partial}{\partial u}} [f\cdot dx_1\cdots dx_N]$ as in the above proof yields
\[  \nabla_{u\frac{\partial}{\partial u}} [f\cdot dx_1\cdots dx_N]= -\frac{\deg( [f\cdot dx_1\cdots dx_N])}{2} \cdot  [f\cdot dx_1\cdots dx_N]. \]
Thus, the homogeneity condition $S3.$ in Definition~\ref{defi-splitting} is indeed to require that the splitting map $s: {\sf Jac} (W) dx_1\cdots dx_N\ra H^*\big( \Omega_{\C^N}^*[[u]], dW+ud_{DR}\big)$ be degree preserving. The existence of a splitting in the rationally graded case is proved in general in~\cite{LLS}. In fact in {\em loc. cit.} the authors wrote down a non-empty parameter space of the set of splittings. The $\omega$-compatibility condition follows from the above Lemma~\ref{lem:omega}. When writing down a homogeneous splitting for $f\cdot dx_1\cdots dx_N$ in general, it looks like
\[ f\cdot dx_1\cdots dx_N+u\cdot  f_1\cdot dx_1\cdots dx_N +u^2 \cdot f_2\cdot dx_1\cdots dx_N +...\]
with $\deg(f_k\cdot dx_1\cdots dx_N)-\deg(f\cdot dx_1\cdots dx_N)=2k$. Assume a lifting of the form $$\sum_{0\leq j\leq k-} f_j dx_1\cdots dx_N \cdot u^j$$ up to $u^{k-1}$ is given, the set of liftings to $u^k$ is a torsor over $\big( {\sf Jac} (W) dx_1\cdots dx_N \big)_{2k+\deg(f\cdot dx_1\cdots dx_N)}$, the degree $\big(\deg(f\cdot dx_1\cdots dx_N)+2k\big)$-part  of the Hochschild homology. Thus to prove Theorem~\ref{thm:canonical} in the case when the group $G=\id$, it suffices to prove the following

\begin{Proposition}~\label{prop:vanishing}
Let $W$ be an invertible polynomial with maximal diagonal symmetry group $G^{{\sf max}}_W$. Then for any degree $a,b\in \Q$ such that $a\neq b$, we have 
\[ {{\sf Hom}}_{\mathbb{K}[G^{{\sf max}}]}\big( ({\sf Jac} (W) dx_1\cdots dx_N )_a, ({\sf Jac} (W) dx_1\cdots dx_N )_b\big)=0.\]
In other words, there exists no non-trivial $\mathbb{K}[G^{{\sf max}}]$-equivariant maps between homogeneous components of different degrees.
\end{Proposition}

\begin{Proof}

It is known that every invertible polynomial $W$ decomposes as
\[ W=W_1+\cdots+W_l,\]
with each $W_l$ so-called atomic polynomials which are of three types: Fermat, Chain, or Loop. We first prove the proposition for each of the three basic types.

{{\bf (i) Fermat type.}} Let $W=x^n$. The weight $q=\frac{1}{n}$. The Hochschild homology is given by
\[ {{\sf Jac}}(W)dx =\langle dx, \cdots, x^{n-2}dx \rangle,\]
with the degree given by $\deg( x^kdx )= \frac{n-2k-2}{n}$ for $0\leq k\leq n-2$. In this case $G^{{\sf max}}_W=\langle g \rangle$ is generated by $g=e^{2\pi i /n}$. And the generator $g$ acts by
\[ g(x^kdx) = e^{2\pi i \cdot \frac{k+1}{n}} x^kdx.\]
Thus $g$ acts on the space of $\mathbb{C}$-linear maps ${{\sf Hom}}\big( ({\sf Jac} (W) dx )_a, ({\sf Jac} (W) dx )_b\big)$ by $e^{2\pi i\cdot \frac{a-b}{2}}$. But since the difference between two different weights $a$ and $b$ satisfies that
\[-1< \frac{a-b}{2}<1,\]
the factor $e^{2\pi i\cdot \frac{a-b}{2}}\neq 1$ for any $a\neq b$. This settles the Fermat case.

{{\bf (ii) Chain type.}} For the Chain and Loop type, the exponent matrix and its inverse plays an important role in the proof. We refer to~\cite[Section 5.1]{HLSW} for details used here. Let $W=x_1^{a_1}+x_1x_2^{a_2}x_3+\cdots + x_{N-1} x_N^{a_N}$ be a chain type invertible polynomial. Denote by $E_W$ its exponent matrix
\[ \begin{bmatrix}
a_1 & 0  &  \cdots & 0 & 0\\
1 &    a_2   & \cdots & 0&0\\
\vdots & \vdots  & \cdots &  \vdots & \vdots\\
0 & 0 & \cdots & 1 & a_N
\end{bmatrix}\]
For a row vector $r=(r_1,\cdots, r_N)$, we write $x^r:= x_1^{r_1}\cdots x_N^{r_N}$. We use a basis of ${{\sf Jac}}(W)dx$ by Krawitz~\cite{Kra} to perform the computation. This special basis consists of monomials $x^rdx$ with $0\leq r_j\leq a_j-1$ and $r\neq (*,\cdots,*,k, a_{N-2l}-1, \cdots, 0, a_{N-2}-1,0, a_N-1)$ with $k\geq 1$. The weights $q_j$'s are related to the matrix $E_W^{-1}$ by $(q_1,\ldots,q_N)^{{\sf t}}= E_W^{-1}  \cdot (1,1,1,\cdots,1,1)^{{\sf t}}$, which implies that
\[ \frac{\deg(x^{r'}dx)-\deg(x^{r}dx)}{2} = (r-r')\cdot E_W^{-1} \cdot (1,1,1,\cdots,1,1)^{{\sf t}}.\]
Furthermore, the $G^{{\sf max}}$-equivariant condition  is equivalent to that $(r-r')\cdot E_W^{-1}\in \Z^N$ since the column vectors of $E_W^{-1}$ generates the group $G^{{\sf max}}$. Without loss of generality, let us assume that $\deg(x^{r'}dx)>\deg(x^{r}dx)$. Then we have that 
\[ \frac{\deg(x^{r'}dx)-\deg(x^{r}dx)}{2} = (r-r')\cdot E_W^{-1} \cdot (1,1,1,\cdots,1,1)^{{\sf t}} \in \Z_{>0}.\]
This implies that at least one entry in $(r-r')\cdot E_W^{-1}$ is greater or equal to $1$. But on the other hand, we may explicitly compute the matrix $E_W^{-1}$:
\[ \big(E_W^{-1}\big)_{ij} = (-1)^{i+j} \prod_{j\leq l \leq i} \frac{1}{a_l}.\]
If the $j$-th entry in $(r-r')\cdot E_W^{-1}$ is greater or equal to $1$, we get
\begin{equation}~\label{eq:geq1}
 \sum_{j\leq i\leq N} (-1)^{i+j}  \frac{(r_i-r_i')}{\prod_{j\leq l \leq i} a_l}  \geq 1.
 \end{equation}
However, since $|r_i-r_{i'}|\leq a_i-1$, we also have
\begin{align*}
|\sum_{j\leq i\leq N} (-1)^{i+j}  \frac{(r_i-r_i')}{\prod_{j\leq l \leq i} a_l} | & \leq \sum_{j\leq i\leq N} \frac{(a_i-1)}{\prod_{j\leq l \leq i} a_l} \\
&= \frac{a_j-1}{a_j}+\frac{a_{j+1}-1}{a_ja_{j+1}}+\frac{a_{j+2}-1}{a_ja_{j+1}a_{j+2}}+\cdots\\
&= \frac{a_ja_{j+1}-1}{a_ja_{j+1}}+\frac{a_{j+2}-1}{a_ja_{j+1}a_{j+2}}+\cdots\\
&=\frac{a_ja_{j+1}a_{j+2}-1}{a_ja_{j+1}a_{j+2}}+\cdots\\
&< 1.
\end{align*}
Thus the inequality~\ref{eq:geq1} cannot hold, which proves the statement for the Chain type atomic polynomials.

{{\bf (iii) Loop type.}} In this case, the exponent matrix is given by
\[ \begin{bmatrix}
a_1 & 0  &  \cdots & 0 & 1\\
1 &    a_2   & \cdots & 0&0\\
\vdots & \vdots  & \cdots &  \vdots & \vdots\\
0 & 0 & \cdots & 1 & a_N
\end{bmatrix}\]
Its determinant is $D=\prod_{k=1}^N a_k-(-1)^N$. Its inverse is 
\[ (E_W^{-1})_{ij}= \begin{cases}
(-1)^{N+i+j} \frac{\prod_{k=i+1}^{j-1}a_k}{D}, & \; i<j\\
(-1)^{i+j} \frac{\prod_{k=i+1}^N a_k \cdot \prod_{l=1}^{j-1} a_l }{D}, & \; i\geq j
\end{cases},\]
with the convention that empty product is $1$. As in the chain type case, we calculate using two basis elements $x^rdx$ and $x^{r'}dx$, with $r=(r_1,\ldots,r_N)$ and $r'=(r'_1,\ldots,r'_N)$ such that $0\leq r_j,r_j'\leq a_j-1, \; (\forall 1\leq j\leq N)$. As in the chain type, it suffices to show that every entry in $(r-r')\cdot E_W^{-1}$ is strictly less than $1$. Indeed, we have its $j$-th entry is
\begin{align*}
& \sum_{1\leq i< j} (-1)^{N+i+j} (r_i-r_i')  \frac{\prod_{k=i+1}^{j-1}a_k}{D} +\sum_{j\leq i\leq N} (-1)^{i+j}(r_i-r_i') \frac{\prod_{k=i+1}^N a_k \cdot \prod_{l=1}^{j-1} a_l }{D} \\
\leq & \frac{ \prod_{k=1}^{j-1}a_k - \prod_{k=2}^{j-1} a_k +  \cdots +a_{j-1}-1}{D} + \frac{ \big( \prod_{k=j}^N a_k - \prod_{k=j+1}^N a_k +\cdots -1\big) \prod_{l=1}^{j-1} a_l}{D}\\
= &  \frac{ (\prod_{k=1}^{j-1}a_k -1 ) + (\prod_{k=j}^N a_k -1)\cdot  \prod_{l=1}^{j-1} a_l}{D}\\
= & \frac{ \prod_{k=1}^N a_k -1 }{D}
\end{align*}
If $N$ is odd, the fraction $\frac{ \prod_{k=1}^N a_k -1 }{D}= \frac{ \prod_{k=1}^N a_k -1 }{\prod_{k=1}^N a_k +1}<1$, and we are done. If $N$ is even, the fraction $\frac{ \prod_{k=1}^N a_k -1 }{D}=1$. In the above inequality, equality can only hold when (with $N$ even)
\[ (-1)^{i+j}(r_i-r_i')= a_i -1, \;\; \forall 1\leq i\leq N.\]
If this is the case, one then verifies that 
\[ (r-r')\cdot E_W^{-1} = ((-1)^{j+1},(-1)^j,\ldots,1,\ldots,(-1)^{j+N}),\]
i.e. the right hand side vector is the unique alternating vector with entries $\pm 1$ and the $j$-th position is $1$. But, then this implies that $ \frac{\deg(x^{r'}dx)-\deg(x^{r}dx)}{2} = (r-r')\cdot E_W^{-1} \cdot (1,1,1,\cdots,1,1)^{{\sf t}}=0$ since $N$ is even, which contradicts the assumption that $\deg(x^rdx)\neq \deg(x^{r'}dx)$. This proves the case for the Loop type invertible polynomails.

To finish the proof, for a general invertible polynomial $W=W_1+\cdots+W_l$, we observe that
\[ G_W^{{\sf max}} = G_{W_1}^{{\sf max}} \times \cdots \times G_{W_l}^{{\sf max}}.\]
Furthermore, the Hochschild homology also decompose into
\[ {{\sf Jac}}(W)={{\sf Jac}}(W_1)\otimes\cdots\otimes {{\sf Jac}}(W_l).\]
Using that $G_W^{{\sf max}}$ is commutative, the result easily follows from the atomic cases.
\end{Proof}

\paragraph{{\bf Canonical splitting with an orbifolding group $G\subset G_W^{{\sf max}}\cap \SL_N$.}} To prove Theorem~\ref{thm:canonical} with an orbifolding group $G$, the idea is to use localization formula of Hochschild invariants to reduce the statement to each invariant pieces. For this purpose, we need the following 

\begin{Lemma}~\label{lem:invariants}
Let $W$ be an invertible polynomial. Let $g\in G^{{\sf max}}$ be a diagonal symmetry of $W$. Denote by $V^g$ the fixed point locus of the action $g: V\ra V$. Then $W|_{V^g}$ is also an invertible polynomial.
\end{Lemma}

\begin{Proof}
Assume that $W=W_1+\cdots+W_l$ decomposes into sums of atomic polynomial $W_j$'s. Since $g$ is a diagonal symmetry, it also decomposes as $g=(g_1,\cdots, g_l)$, which shows that the restriction $W|_{V^g}$ is the sum of the restrictions $W_j$ on the invariant subspace of $g_j$. Thus it suffices to prove the statement for atomic polynomials. 
\begin{itemize}
\item For the Fermat type $W=x^n$, the fixed subspace is either $\C$ or $0$, thus clearly invertible. 
\item For the chain type $W=x_1^{a_1}+x_1x_2^{a_2}+\cdots + x_{N-1} x_N^{a_N}$, notice that if $g(x_i)=x_i$ with $i\geq 2$, then it implis that $g(x_{i-1})=x_{i-1}$. Thus if we set $l(g)$ to be the maximal index such that $g(x_{l(g)})=x_{l(g)}$, we have that
\[ V^g=\C^{l(g)}_{x_1,\cdots,x_{l(g)}}.\]
It is clear that $W|_{V^g}$ is again atomic and of chain type. 
\item For the loop type, one can show that the fixed subspace $V^g$ is either the whole space $V$ or $0$, thus clearly invertible. 
\end{itemize}
\end{Proof}

\paragraph{{\bf Finishing the proof of Theorem~\ref{thm:canonical}.}} Recall the following localization formula:
\begin{equation}~\label{eq:localization}
 HH_*(\MF_G(W))\cong \bigoplus_{g\in G} \big( HH_*(\MF(W|_{V^g}))\big)_{G}.
 \end{equation}
Similarly for its negative cyclic homology. Again, since $G_W^{\sf max}$ is a diagonal symmetry group, requiring a splitting to be $G_W^{\sf max}$-equivariant is the same as requiring the splitting be $G_{W|_{V^g}}^{{\sf max}}$-equivariant on each $g$-invariant subspace. Now the uniqueness follows from that on each $g$-invariant subspace by Lemma~\ref{lem:invariants}.

\section{The cubic family}~\label{sec:cubic}

As an example, we work out the case of the Fermat cubic
\[ W=\frac{1}{3}(x_1^3+x_2^3+x_3^3)\]
endowed with the orbifold group $G=\Z/3\Z=\langle (\zeta_3,\zeta_3,\zeta_3)\rangle$ acting diagonally on the $x$-variables, with $\zeta_3=e^{2\pi i/3}$. The purpose of the section is to demonstrate the computability of non-commutative Hodge structures and categorical primitive forms in the simplest Calabi-Yau case. Overall, we have chosen to provide a schematic flow of how the non-commutative computation could be done, while tedious details of most computations are omitted. 

To write down explicit formulas, it is convenient to use the Shuffle and cyclic Shuffle operations, denoted by $\shu$ and $\shu^c$ respectively. We refer to~\cite[Section 2.2]{Shk2} for the definitions of them.

\paragraph{{\bf The $A_\infty$-algebra structure.}} By our previous work~\cite{Tu}, by Kontsevich's deformation quantization formula one can write down an $A_\infty$-algebra structure on the $Z/2\Z$-graded vector space
\[ {{\sf End}}(\C^{\sf stab})= \Lambda^*(\epsilon_1,\epsilon_2,\epsilon_3),\]
i.e. exterior tensors generated by three odd variables $\epsilon_1, \epsilon_2, \epsilon_3$. We denote by $\epsilon_{12}=\epsilon_1\wedge\epsilon_2$, and similarly $\epsilon_{213}=\epsilon_2\wedge\epsilon_1\wedge \epsilon_3$. We shall denote this $A_\infty$-algebra by $E$.

For the Hochschild invariants of $E$, one can use the inverse of the HKR-map to obtain the following

\begin{Proposition}~\label{prop:basis}
The Hochschild homology of $E$ is an $8$-dimensional vector space, concentrated at the odd degree. Explicitly, a basis is given by
\begin{align*}
HH_*(E) & = {{\sf span}}\big( [ \epsilon_{123}], [\epsilon_{123}|\epsilon_1],[\epsilon_{123}|\epsilon_2], [\epsilon_{123}|\epsilon_3],\\ &[\epsilon_{123}|\shu (\epsilon_1,\epsilon_2)], [\epsilon_{123}|\shu(\epsilon_1,\epsilon_3)], [ \epsilon_{123}|\shu(\epsilon_2,\epsilon_3)], [\epsilon_{123}|\shu(\epsilon_1,\epsilon_2,\epsilon_3)]\big).
\end{align*}
\end{Proposition} 

\medskip
\begin{remark}
The above formula of generators of the Hochschild homology can also be obtained using the Kunneth product~\cite{Shk}. Indeed, let $A$ and $B$ be two strictly unital $A_\infty$ algebras. We have the Kunneth map $\times : C_*(A)\otimes C_*(B) \ra C_*(A\otimes B)$ defined by
\[ (a_0|a_1|\cdots|a_n)\times (b_0|b_1|\cdots|b_m) = (-1)^* a_0\otimes b_0| {{\sf sh}}\big( a_1\otimes\bone|\cdots|a_n\otimes \bone |\bone\otimes b_1|\cdots|\bone\otimes b_m\big).\]
Here the sign $*=|b_0|\cdot (|a_1|'+\cdots+|a_n|')$. One can apply the formula by considering $E\cong (A_3\otimes A_3) \otimes A_3$ with $A_3$ the $A_\infty$ algebra associated with $\frac{1}{3}x^3$. This approach is a bit {\sl ad. hoc.} since the tensor product of $A_\infty$-algebras are rather complicated to define, and may not be an associative operation on triple tensor products. Nevertheless, the Kunneth formula seems to still hold.
\end{remark}

\paragraph{{\bf The semi-direct product $E\rtimes \Z/3\Z$.}} The orbifold group $\Z/3\Z$ acts on $E$ diagonally on the variables $\epsilon_1$, $\epsilon_2$ and $\epsilon_3$. Since Kontsevich's deformation quantization is invariant under the affine transformation group, the $A_\infty$-algebra structure on $E$ is therefor $\Z/3\Z$-equivariant. Thus, we may form the semi-direct product $A_\infty$-algebra $E\rtimes \Z/3\Z$. By general orbifolding construction of compact generators, this latter algebra is Morita equivalent to the category $\MF_G(W)$.

The Hochschild homology of $E\rtimes \Z/3\Z$ can be computed using the localization formula~\ref{eq:localization}. In this formula, we used the $G$-coinvariants which is isomorphic to the $G$-invariants since our base field $\mathbb{K}$ has characteristic zero. Throughout the section, we shall work with $G$-invariants instead of coinvariants. Since the dimension $N=3$ is odd, Formula~\ref{eq:localization} implies that
\[ HH_{\sf odd}(E\rtimes \Z/3\Z) \cong \big( HH_{\sf odd}(E) \big)^{\Z/3\Z} ={{\sf span}} \big( [ \epsilon_{123}],[\epsilon_{123}|\shu(\epsilon_1,\epsilon_2,\epsilon_3)]\big).\]

\paragraph{{\bf Explicit formula of the canonical splitting.}}  The maximal diagonal symmetry group of $W$ is given by $G_W^{\sf max}=(\Z/3\Z)^3$. We may write down an explicit formula for the unique $G_W^{{\sf max}}$-equivariant splitting of $E\rtimes \Z/3\Z$. These formulas are obtained through the so-called cyclic Kunneth map~\cite{Shk2}.  Explicitly, the canonical splitting is given by
\begin{align}~\label{eq:s}
\begin{split}
s([\epsilon_{123}])&=\sum_{i\geq 0, j\geq 0, k\geq 0} (-1)^{i+j+k} d_id_jd_k\; \epsilon_{123}|\shu ( \epsilon_1^{3i}, \epsilon_2^{3j}, \epsilon_3^{3k}) u^{i+j+k}\\
&+\sum_{i\geq 0, j\geq 0, k\geq 0} (-1)^{i+j+k+1}d_id_jd_k\; \epsilon_3|\shu \big( \shu^c(\epsilon_1^{3i+1},\epsilon_2^{3j+1}), \epsilon_3^{3k}\big) u^{i+j+k+1}\\
&+ \sum_{i\geq 0, j\geq 0, k\geq 0} (-1)^{i+j+k}d_id_jd_k \; \bone|\shu^c\big( \epsilon_{12}|\shu(\epsilon_1^{3i}, \epsilon_2^{3j}), \epsilon_3^{3k+1}\big) u^{i+j+k+1}.
\end{split}
\end{align}
Here $d_i:=\prod_{1\leq l\leq i} (3l-2)$, and $d_0=1$. Another similar formula looks like
\begin{align}~\label{eq:omega}
\begin{split}
& s([\epsilon_{123}|\shu(\epsilon_1,\epsilon_2,\epsilon_3)])\\
=&\sum_{i\geq 0, j\geq 0, k\geq 0} (-1)^{i+j+k} c_ic_jc_k\; \epsilon_{123}|\shu ( \epsilon_1^{3i+1}, \epsilon_2^{3j+1}, \epsilon_3^{3k+1}) u^{i+j+k}\\
&+\sum_{i\geq 0, j\geq 0, k\geq 0} (-1)^{i+j+k+1}c_ic_jc_k\; \epsilon_3|\shu \big( \shu^c(\epsilon_1^{3i+2},\epsilon_2^{3j+2}), \epsilon_3^{3k+1}\big) u^{i+j+k+1}\\
&+ \sum_{i\geq 0, j\geq 0, k\geq 0} (-1)^{i+j+k}c_ic_jc_k \; \bone|\shu^c\big( \epsilon_{12}|\shu(\epsilon_1^{3i+1}, \epsilon_2^{3j+1}), \epsilon_3^{3k+2}\big) u^{i+j+k+1},
\end{split}
\end{align}
with $c_i:=\prod_{1\leq l\leq i} (3l-1)$, and $c_0=1$. One can check directly that the map $s$ defined above is $G_W^{\sf max}$-equivariant with respect to the diagonal action of $G_W^{\sf max}$ on the elements $\epsilon_1$, $\epsilon_2$ and $\epsilon_3$. Hence by Theorem~\ref{thm:canonical}, this is the canonical splitting of the non-commutative Hodge filtration.

\paragraph{{\bf Deformations of the Fermat cubic.}} By the localization formula for the Hochschild cohomology, we have
\[ HH^{{\sf even}}(E)^{\Z/3\Z}\cong HH^{{\sf even}} (E\rtimes \Z/3\Z).\]
Hence all deformations of $E\rtimes \Z/3\Z$ are from $\Z/3\Z$-invariant deformations of $E$. A universal family can be realized by deforming $W$ to
\[ \mathcal{W}:= W+t_0+t\cdot x_1x_2x_3,\]
and then apply Kontsevich's deformation quantization to $\mathcal{W}$ to obtain a family of $A_\infty$-algebras $E_\mathcal{W}$. From this family of $A_\infty$-algebras, taking negative cyclic homology yields a VSHS on $HC_*^-(E_\mathcal{W})$ over the formal power series ring $\C[[t_0,t]]$. One can show that there is an isomorphism of VSHS's~\footnote{See Proposition~\ref{prop:equiv} for an analogous statement and its proof.}:
\[ \cV^{\MF_{\Z/3\Z}(W)} \cong HC_*^-(E_\mathcal{W})^{\Z/3\Z}.\]
In the remaining part of the section, we use the above isomorphism and Theorem~\ref{thm:bijection2} to compute the primitive form $\zeta$ associated with the canonical splitting in Equations~\ref{eq:s} and~\ref{eq:omega}. Since the $t_0$-direction does not affect the main computation, we shall set $t_0=0$ in the remaining part of the section. We begin with the Kodaira-Spencer class of the above family. Denote by
\[ \mu(t)_n: (E_\mathcal{W}[1])^{\otimes n} \ra E_\mathcal{W}[1]\]
the $A_\infty$-structure obtained from deformation quantization. Then we have
\begin{Lemma}
For any multi-index $J$, we have the following identities:
\begin{align*}
\frac{\partial}{\partial t}\mu(t)_3\big(\epsilon_{\sigma(1)}\wedge \epsilon_J,\epsilon_{\sigma(2)},\epsilon_{\sigma(3)}\big)&= \frac{1}{6}\epsilon_J, \;\;\forall \sigma\in \Sigma_3.\\
\frac{\partial}{\partial t}\mu(t)_n\big(\epsilon_{i_1}\wedge \epsilon_J,\epsilon_{i_2},\cdots,\epsilon_{i_n}\big)&=0, \;\;\forall n\geq 4, \;\;\forall 1\leq i_j\leq 3.
\end{align*}
\end{Lemma}

The Getzler-Gauss-Manin connection acts on cyclic chains of $E_{\mathcal{W}}$ by formula
\begin{align*}
 \nabla^{\sf Get}_{\frac{\partial}{\partial t}} (\xi)&= \frac{\partial}{\partial t}(\xi)-u^{-1}\iota\big(\frac{\partial\mu(t)}{\partial t}\big)(\xi),\\
\iota\big(\frac{\partial\mu(t)}{\partial t}\big)(\xi)&= b^{1|1}(\frac{\partial\mu(t)}{\partial t}; \xi)+B^{1|1}(\frac{\partial\mu(t)}{\partial t}; \xi).
\end{align*}

For simplicity, we set $A:=\iota(\frac{\partial}{\partial t}\mu(t))$.

\begin{Lemma}
There exists an isomorphism of differential modules over $\mathbb{K}[[t]]$:
\[ e^{tA/u}: \big(HP_*(E)\widehat{\otimes}_\mathbb{K}\mathbb{K}[[t]], \frac{\partial}{\partial t}\big) \ra \big( HP_*(E_{\mathcal{W}}), \nabla^{\sf Get}_{\frac{\partial}{\partial t}}\big).\]
In other words, the above isomorphism trivializes the differential module $HP_*(E_\mathcal{W})$.
\end{Lemma}

\begin{Lemma}~\label{lem:A-action}
For any positive integers $I, J, K$, we have
\begin{align*}
A\big( \epsilon_{123}|\shu(\epsilon_1^I, \epsilon_2^J, \epsilon_3^K)\big)&= \epsilon_{123}|\shu(\epsilon_1^{I-1}, \epsilon_2^{J-1}, \epsilon_3^{K-1}),\\
A\big( \epsilon_3|\shu(\epsilon_1|\epsilon_2, \epsilon_1^I, \epsilon_2^J, \epsilon_3^K)\big) &=(I+\frac{1}{2})\cdot\epsilon_3|\shu(\epsilon_1^I,\epsilon_2^J,\epsilon_3^{K-1})+\epsilon_3|\shu(\epsilon_1|\epsilon_2,\epsilon_1^{I-1}, \epsilon_2^{J-1}, \epsilon_3^{K-1})\\
A\big( \bone|\shu(\epsilon_{12}|\epsilon_3, \epsilon_1^I, \epsilon_2^J, \epsilon_3^K)\big)&= \frac{1}{2}\cdot \epsilon_1|\shu(\epsilon_1^{I-1},\epsilon_2^J,\epsilon_3^K)-\frac{1}{2}\cdot \epsilon_2|\shu(\epsilon_1^I,\epsilon_2^{J-1},\epsilon_3^K)\\
&+  \bone |\shu(\epsilon_{12}|\epsilon_3, \epsilon_1^{I-1},\epsilon_2^{J-1},\epsilon_3^{K-1})+\bone |\shu (\epsilon_{12}, \epsilon_1^{I-1}, \epsilon_2^{J-1}, \epsilon_3^K).
\end{align*}
\end{Lemma}

\paragraph{{\bf Computation of the primitive form.}} Let us write
\[ s_0:= s\big(\epsilon_{123}\big),\;\; \omega:= s\big(\epsilon_{123}|\shu (\epsilon_1,\epsilon_2,\epsilon_3)\big)\] 
for the negative cyclic homology classes of $E$ in Equation~\ref{eq:s} and~\ref{eq:omega}. The class $\omega$ induces a trace map
\[ \Tr_\omega: E \ra \mathbb{K}, \;\; \Tr_\omega(\epsilon_{123})=1, \mbox{\; and $0$ otherwise.}.\]
The graded symmetric bilinear form $\langle\epsilon_I,\epsilon_J\rangle=\Tr_\omega(\epsilon_{I}\wedge\epsilon_J)$ is non-degenerate, which implies that $\omega$ defines a   Calabi-Yau structure of $E$. As we have proved in the previous section, the canonical splitting is $\omega$-compatible. 

The canonical splitting $s$ determines a primitive form of the VSHS on $HC_{{\sf odd}}^-(E_{\mathcal{W}})^{\Z/3\Z}$ by Theorem~\ref{thm:bijection2}. We briefly recall this construction. Indeed, the splitting $s$ defines a subspace
\[ L^s:= \bigoplus_{l\geq 1} u^{-l} \cdot \imag (s) \subset \big(HP_*(E)\big)^{\Z/3\Z},\]
complimentary to the canonical subspace $(HC_*^-(E))^{\Z/3\Z}\subset (HP_*(E))^{\Z/3\Z}$. We can parallel transport the subspace $L^s$ using the Getzler-Gauss-Manin connection to nearby fibers, which yields a direct sum decomposition 
\[ HP_*(E_\mathcal{W})^{\Z/3\Z} \cong HC_*^-(E_\mathcal{W})^{\Z/3\Z}\oplus (L^s)^{{\sf flat}}.\]
Here and in the following, we use an superscript ${{\sf flat}}$ to mean Getzler-Gauss-Manin flat extension of periodic cyclic homology classes. Equivalently, this direct sum decomposition corresponds to a projection map 
\[ \Pi^s : HP_*(E_\mathcal{W})^{\Z/3\Z} \ra HC_*^-(E_\mathcal{W})^{\Z/3\Z}\]
which splits the canonical inclusion map. Then, the primitive form associated with the splitting $s$ is given by 
\[ \zeta:= \Pi^s\big( \omega^{{\sf flat}}\big).\] 
That is, the primitive form $\zeta$ is simply the projection of the flat extension $\omega^{{\sf flat}}$ which {{\sl a priori}} is only a periodic cyclic homology class, along the parallel transported splitting.

Using Lemma~\ref{lem:A-action} to compute $s_0^{{\sf flat}}$ and $\omega^{{\sf flat}}$ yields
\begin{align*}
s_0^{{\sf flat}} &= \sum_{n\geq 0} \frac{ (-1)^n (d_n)^3 t^{3n}}{(3n)!} \epsilon_{123} + O(u),\;\; d_n=\prod_{1\leq i\leq n} (3i-2), \;\; d_0=1\\
\omega^{{\sf flat}} &= \sum_{n\geq 0} \frac{ (-1)^n (c_n)^3 t^{3n+1}}{(3n+1)!} \epsilon_{123} u^{-1} + O(u^0),\;\; c_n=\prod_{1\leq i\leq n} (3i-1), \;\; c_0=1
\end{align*}
Let us set $g(t):=\sum_{n\geq 0} \frac{ (-1)^n (d_n)^3 t^{3n}}{(3n)!}$ and $h(t):=\sum_{n\geq 0} \frac{ (-1)^n (c_n)^3 t^{3n+1}}{(3n+1)!}$. Observe that $g(t)=1+O(t)$, hence it is invertible in $\mathbb{K}[[t]]$. Using $s_0^{{\sf flat}}$ to kill the $u^{-1}$-term of $\omega^{{\sf flat}}$ shows that
\[ \zeta:= \omega^{{\sf flat}} - \frac{h(t)}{g(t)} s_0^{{\sf flat}} u^{-1} \in O(u^0).\]

\paragraph{{\bf The flat coordinate $\tau$.}} The primitive form induces a flat structure on the deformation space. In this flat structure, as pointed out in~\cite[Section 3.2.2]{LLSS}, the flat coordinate can be read off from the coefficient of $u^{-1}$-term in $\omega^{{\sf flat}}$ which gives $\tau= \frac{h(t)}{g(t)}$. Using this flat coordinate to compute the genus zero prepotential function by formula
\[ \partial_i\partial_j\partial_k \cF := \langle  u\nabla_{\partial_i} u \nabla_{\partial_j} u \nabla_{\partial_k} \zeta, \zeta \rangle_{\hres},\]
with the vector fields $\partial_i$, $\partial_j$, $\partial_k$ are either $\partial/\partial t_0$ or $\partial/\partial \tau$. We obtain that  $\cF= \frac{1}{2}t_0^2\tau$, confirming that for the category $\MF_G(W)$ (which, by Orlov~\cite{Orl}, is equivalent to the derived category of coherent sheaves on the Elliptic curve $W=0$ in $\CP^2$) there is no higher non-vanishing invariants in genus zero.

\section{The quintic family}~\label{sec:quintic}

In this section, we study VSHS's associated with quintic families. A notable difference with the discussion in Section~\ref{sec:hodge} is that in this section, we shall work with VSHS's over an actual geometric base $B$, rather than a formal base. Another difference is that we no longer require primitivity. We shall emphasize these differences whenever it is necessary. 

Recall the mirror quintic family $[\mathfrak{X}/(\Z/5\Z)^3]$ is defined as follows.  First, consider a one-parameter family of smooth Calabi-Yau hypersurfaces $\mathfrak{X}\subset B \times \CP^4$ defined by the following equation
\[ \mathcal{W}(\psi,x):= \frac{1}{5}( x_1^5+x_2^5+x_3^5+x_4^5+x_5^5)-\psi x_1x_2x_3x_4x_5,\]
with $\psi\in B=\C\backslash\{1,e^{2\pi i/5}, e^{4\pi i/5}, e^{6\pi i/5}, e^{8\pi i/5}\}$. Denote by $\pi: \mathfrak{X}\ra B$ the canonical the projection map. Next, we observe that the group 
\[ (\Z/5\Z)^3=\{ (a_1,a_2,a_3,a_4,a_5)\in (\Z/5\Z)^5 \mid a_1+\cdots+a_5=0, a_5=0\}\]
naturally acts on $\mathfrak{X}$ by sending $x_j \; (j=1,\ldots,5)$ to $e^{2 a_j \pi i/5} \cdot x_j$. We thus obtain a family of smooth Calabi-Yau orbifolds $[\mathfrak{X}/(\Z/5\Z)^3]$. 

Recall the following Lemma from~\cite[Lemma 2.7]{GPS}.

\begin{Lemma}~\label{lem:equivalence}
A $\Z$-graded polarized VSHS $(\cE,\nabla,\langle-,-\rangle_{{\sf hres}})$ of parity $N$ over a base space $B$ is equivalent to the following data:
\begin{itemize}
\item[--] A locally free, finite rank, $\Z/2\Z$-graded $S$-module $\cV\cong \cV_{[0]}\oplus \cV_{[1]}$.
\item[--] A flat connection $\nabla$ on each $\cV_\sigma$.
\item[--] A decreasing filtration $F^\bullet \cV_{\sigma}$ on each $\cV_\sigma$ that satisfies the Griffiths' transversality $\nabla F^p\cV_{\sigma}\subset F^{p-1}\cV_{\sigma}$.
\item[--] A covariantly constant bilinear pairing
\[ (-,-): \cV_\sigma\otimes \cV_\sigma \ra S\]
such that $(\alpha,\beta)=(-1)^N(\beta,\alpha)$, and satisfies the property that $(F^p\cV_\sigma,F^q\cV_\sigma)=0$ if $p+q>0$, and the induced pairing $(-,-): {\sf Gr}^p_F\cV_\sigma\otimes {\sf Gr}^{-p}_F\cV_\sigma \ra S$ is non-degenerate for all $p$.
\end{itemize}
\end{Lemma}

\begin{proof}
We refer to~\cite[Lemma 2.7]{GPS} for the proof. For later use, we recall the following from {\sl loc. cit.}. First, we form the localization bundle $\widetilde{\cE}:=\cE[u^{-1}]$. Notice that we have the periodicity isomorphism $\widetilde{\cE}_k \stackrel{u\cdot}{\longrightarrow} \widetilde{\cE}_{k-2}$.  Associated with the polarized VSHS $(\cE,\nabla,\langle-,-\rangle_{{\sf hres}})$, one sets
\begin{itemize}
\item $\cV_\sigma := \widetilde{\cE}_k$ with any $k$ such that $k\pmod{2}=\sigma$.
\item The bilinear pairing is given by
\begin{equation}~\label{eq:bilinear}
(\alpha,\beta):= {\sqrt{-1}}^{\,k} \langle \tilde{\alpha},\tilde{\beta} \rangle_{\sf hres}
\end{equation}
if $\tilde{\alpha}\in \widetilde{\cE}_k$ and $\tilde{\beta}\in \widetilde{\cE}_{-k}$. Note that the right hand side lies inside $S$ by degree reason.
\end{itemize}
\end{proof}

We shall freely use this equivalence in the rest of the section.  For a $\Z$-graded VSHS $(\cE,\nabla,\langle-,-\rangle_{{\sf hres}})$, to simplify the notations, denote by $\big(\cE_{\sigma},\nabla,(-,-)\big) \; (\sigma\in \{ {[0]}, {[1]}\})$ the corresponding $\Z/2\Z$-graded module obtained via the lemma above.

\paragraph{{\bf Griffiths' construction of VSHS's.}} Back to the orbifold $[\mathfrak{X}/(\Z/5\Z)^3]$, there are two VSHS's naturally associated with this data: one geometric and the other one categorical. We first recall the geometric construction of VSHS, essentially due to Griffiths~\cite{Gri1}. Roughly speaking, one first constructs a resolution~\cite{Roa}, say $\widetilde{\mathfrak{X}}$ of the orbifold $[\mathfrak{X}/(\Z/5\Z)^3]$. Denote by $\widetilde{\pi}: \widetilde{\mathfrak{X}} \ra B$ the projection map to the base. Define a VSHS using the equivalent data  described in Lemma~\ref{lem:equivalence} as follows.  The underlying $\Z/2\Z$-graded bundle is purely odd given by $R^3\widetilde{\pi}_*\C\otimes_\C S$, the middle cohomology of the family $\widetilde{\pi}: \widetilde{\mathfrak{X}}\ra B$. It is endowed with the classical Hodge filtration $F^\bullet$, and the Gauss-Manin connection $\nabla$.  The polarization $(-,-)$ is the intersection form on middle cohomology. It turns out in this case, the VSHS can be described in terms of $\mathfrak{X}$ and taking $(\Z/5\Z)^3$-invariants. Indeed, we have a natural isomorphism 
\begin{equation}~\label{eq:iso1}
 \big( R^3\pi_*\C\otimes_\C S \big)^{(\Z/5\Z)^3} \cong R^3\widetilde{\pi}_*\C\otimes_\C S
 \end{equation}
of VSHS's over $B$. Indeed, consider the pull-back map via $p: \widetilde{\mathfrak{X}}\ra [ \mathfrak{X}/(\Z/5\Z)^3]$. The morphism $p^*$ intertwines the Hodge filtration and the Gauss-Manin connection. For the polarization, we have
\[ (p^*x, p^*y)= 5^3\cdot (x,y), \;\; \forall x, y \in  \big( R^3\pi_*\C\otimes_\C S \big)^{(\Z/5\Z)^3}.\]
The factor $5^3$ is a constant, which shows that $5^{3/2}\cdot p^*$ is an isomorphism between VSHS's.

\paragraph{{\bf Categorical construction of VSHS's.}} The second construction of VSHS over $B$ is analogous to the one described in Section~\ref{sec:hodge}, with the difference that we work over $S=\Gamma(B,\cO_B)$, the ring of functions on $B$, instead of over a formal ring.  More precisley, denote by $A:=\mathbb{C}[\epsilon_1,\epsilon_2,\epsilon_3,\epsilon_4,\epsilon_5]$ the super-commutative ring generated by $5$ odd variables.  Recall as in~\cite{Tu}, we may use Kontsevich's deformation quantization formula to obtain a $S$-linear $A_\infty$-algebra structure on $A\otimes_\mathbb{C} R$. Indeed, let $\cU: T_{\sf poly}(A)[1] \ra  C^*(A)[1]$ be Kontsevich's $L_\infty$ quasi-isomorphism. Tensoring it with $S$ yields (which we still denote by $\cU$)
\[ \cU: T_{\sf poly}(A)[1]\otimes_{\mathbb{C}} S \ra C^*(A\otimes_\mathbb{C} S)[1].\]
The element $\mathcal{W}=\frac{1}{5}( x_1^5+x_2^5+x_3^5+x_4^5+x_5^5)-\psi x_1x_2x_3x_4x_5\in S[x_1,x_2,x_3,x_4,x_5]$ is obviously a Maurer-Cartan element of the left hand side, its push-forward
\[ \cU_*\mathcal{W}:= \sum_{k\geq 1} \frac{1}{k!}\cU_k(\mathcal{W}^k)\]
is a Maurer-Cartan element of $C^*(A\otimes_\mathbb{C} S)[1]$, which by definition is an ($S$-linear) $A_\infty$-algebra structure on $A\otimes_\mathbb{C} S$. We remark that the above infinite sum is well-defined, since for each element $a_1\otimes\cdots\otimes a_n\in A[1]^{\otimes n}$ the evaluation $ \sum_{k\geq 1} \frac{1}{k!}\cU_k(\mathcal{W}^k)(a_1,\ldots,a_n)$ is only a finite sum, since by Kontsevich's explicit formula, each $\mathcal{W}$ acts on the $a$'s by a $5$-th order differential operator. Denote this $S$-linear $A_\infty$-algebra by $A^{\mathcal{W}}$, indicating that it is a deformation of $A$ using the Maurer-Cartan element $\mathcal{W}$. Note that as a $\Z/2\Z$-graded $S$-module, it is simply $A\otimes_\mathbb{C} S$. Since Kontsevich's deformation quantization formula is invariant under any affine symmetry, the maximal diagonal symmetry group $(\Z/5\Z)^4:=\{ (a_1,a_2,a_3,a_4,a_5)\in (\Z/5\Z)^5 \mid a_1+\cdots+a_5=0\}$ acts on the $A_\infty$-algebra $A^\mathcal{W}$. Denote by $A^\mathcal{W}\rtimes (\Z/5\Z)^4$ the associated semi-direct product $A_\infty$-algebra. The relationship of these $A_\infty$-algebras with matrix factorizations is the following. For each fixed $\psi\in B$, it was proved in~\cite{Tu} that the $A_\infty$-algebra $A^{\mathcal{W}_\psi}$ is a minimal model of the dg-algebra ${{\sf End}}(\mathbb{C}^{\sf stab})$ for the Koszul matrix factorization $\mathbb{C}^{\sf stab}\in \MF(\mathcal{W}_\psi)$ which generates the category $\MF(\mathcal{W}_\psi)$ by a result of Dyckerhoff~\cite{Dyc}.

The second VSHS over $B$ is given by the negative cyclic homology $HC_{{\sf odd}}^-\big(  A^\mathcal{W}\rtimes (\Z/5\Z)^4 \big)$. As in Section~\ref{sec:hodge}, we put on it the categorical higher residue pairing, the $u$-connection and the Getzler connection. Similar to the geometric case, one can realize this VSHS as the  $(\Z/4\Z)^4$-invariant part of the VSHS $HC^-(A^\mathcal{W})$.

\begin{Proposition}~\label{prop:equiv}
There exists an isomorphism of VSHS's over $B$:
\begin{equation}~\label{eq:iso2}
 HC_{\sf odd}^-\big( A^\mathcal{W}\rtimes (\Z/5\Z)^4 \big)\cong \big( HC^-( A^\mathcal{W})\big)^{(\Z/5\Z)^4}.
 \end{equation}
 \end{Proposition}

 \begin{proof}
For simplicity, denote $G=(\Z/5\Z)^4$, and its dual group $G^\vee={\sf Hom}(G,\mathbb{C}^*)$ its dual group. There exists a quasi-isomorphism of chain complexes:
 \begin{align*}
\gamma: \big(C_*(A^\mathcal{W}\rtimes G)\big)_{G^\vee} & \ra C_*(A^\mathcal{W})_G\\
 \gamma ( a_0\rtimes g_0|\cdots|a_n\rtimes g_n) &= 
 \begin{cases} 
 a_0|g_0(a_1)|\cdots|g_0g_1\ldots g_{n-1}(a_n) \mbox{\;\; if $g_0\cdots g_n=\id$},\\
 0\mbox{\;\; otherwise}
 \end{cases}
 \end{align*}
The fact $\gamma$ is a quasi-isomorphism follows from the localization formula~\ref{eq:localization} by observing that the condition $g_0\cdots g_n=\id$ is precisely to pick up the component $g=\id$ in this formula. It is straightforward to check that $\gamma$ intertwines both the $u$-connection and the Getzler connection in the $\psi$-direction. For the pairing, we verifies that
\[ |G|\cdot  \langle \gamma(x), \gamma(y) \rangle_{\sf hres}=  \langle x, y \rangle_{\sf hres}.\]
Thus $|G|^{1/2}\cdot \gamma$ is an isomorphism of VSHS's.
 \end{proof}
 
The VSHS $HC^-( A^\mathcal{W})$ can be compared with Saito's original VSHS on the twisted de Rham complex. The following result was proved in our previous work~\cite{Tu} in the formal setting. We reprove it here over the ring $S$.

 \begin{Lemma}~\label{lem:comparison}
 There exists an isomorphism of VSHS's:
 \[ \Psi: HC^-(A^\mathcal{W}) \cong H^*\big( \Omega_{B\times \C^5/B}^*[[u]], d\mathcal{W}+ud_{DR}\big)\]
of VSHS's, where the right hand side is Saito's VSHS with connection operators given by $\nabla_{\frac{d}{d\psi}}:=  \frac{d}{d\psi} - \frac{1}{u} \cdot x_1x_2x_3x_4x_5$, and $\nabla_{\frac{d}{du}}:= \frac{d}{du}-\frac{5}{2u}-\frac{\mathcal{W}}{u^2}.$
\end{Lemma}

\begin{proof}
For simplicity, denote by $\cV^{{\sf Saito}}=H^*\big( \Omega_{B\times \C^5/B}^*[[u]], d\mathcal{W}+ud_{DR}\big)$. Reall from~\cite{Tu} the isomoprhism $\Psi$ is, up to scalar, given by the following composition
\[ HC^-(A^\mathcal{W}) \stackrel{\cU^{{\sf Sh},\mathcal{W}}_0}{\longrightarrow} H^*\big( \Omega_A^*\otimes_{\mathbb{C}} R[[u]], L_{\mathcal{W}}+ud_{DR}\big) \ra (\cV^{{\sf Saito}})^\vee \leftarrow \cV^{{\sf Saito}}.\]
The second arrow is always an isomorphism. The third arrow is induced from Saito's higher residue pairing, which is an isomorphism due to its non-degeneracy. The first map $\cU^{{\sf Sh},\mathcal{W}}_0$ is given by a deformation of the Tsygan formality map using the Maurer-Cartan element $\mathcal{W}$. It is standard in deformation theory to show that it is an isomorphism in the formal setting. However, in our setting, we argue a stronger result that it remains an isomorphism over $S[[u]]$. Indeed, first observe that the map
\[ \cU^{{\sf Sh},\mathcal{W}}_0=\sum_{k\geq 0} \frac{1}{k!} \cU^{\sf Sh}_k(\mathcal{W}^k)\]
is in fact a finite sum when applied to any element $a_0|a_1|\cdots|a_n\in A\otimes A[1]^{\otimes n}$, because $\mathcal{W}$ acts by a $5$-order differential operator. To argue that $$\cU^{{\sf Sh},\mathcal{W}}_0: C_*(A^\mathcal{W})[[u]] \ra\big( \Omega_A^*\otimes_{\mathbb{C}} R[[u]], L_{\mathcal{W}}+ud_{DR}\big) $$ is a quasi-isomorphism, it suffices to prove that
\[ \cU^{{\sf Sh},\mathcal{W}}_0: C_*(A^\mathcal{W}) \ra \big( \Omega_A^*\otimes_{\mathbb{C}} R, L_{\mathcal{W}}\big)\]
is a quasi-isomorphism of mixed complexes. To this end, we consider an exhaustive filtration on both sides by number of $\epsilon$'s. By definition, both differentials and the map $\cU^{{\sf Sh},\mathcal{W}}_0$ preserves this filtration, hence it induces a map of the associated spectral sequences on both sides. The $E^1$-page of the left hand side can be computed using the Hochschild-Konstant-Rosenberg isomorphism for the algebra $A$, which yields exactly $\big( \Omega_A^*\otimes_{\mathbb{C}} R, L_{\mathcal{W}}\big)$. The $E^1$-page of the right hand side remains itself since $L_\mathcal{W}$ always reduces the number of $\epsilon$'s. Thus, we see that $\cU^{{\sf Sh},\mathcal{W}}_0$ induces an isomorphism at the $E^1$-page, which proves that it is a quasi-isomorhism.
\end{proof}

\begin{Lemma}~\label{lem:Z-grading}
The $\Z/5\Z$-invariant part $HC^-(A^\mathcal{W})^{\Z/5\Z}$ is a $\Z$-graded VSHS.
\end{Lemma}

\begin{proof}
For this, we use the previous comparison result. Using the homogeneous property of $\mathcal{W}$ and calculate as in Lemma~\ref{lem:omega}, we have
\[ \nabla_{u\frac{d}{du}} (\alpha)= -\frac{\deg(\alpha)}{2} \alpha,\]
with the degree calculated by
\begin{equation}~\label{eq:degree}
\deg(x_j) = -2/5,\; \deg(dx_j) = 3/5,\; \deg (u) = -2, \; \deg(\psi) =0
\end{equation}
Now, consider the diagonal action of $\Z/5\Z$ on $\mathcal{W}$.  Since $\mathcal{W}$ has an isolated singularity at $x=0$, the cohomology group $H^*\big( \Omega_{B\times \C^5/B}^*[[u]], d\mathcal{W}+ud_{DR}\big)$ is concentrated at differential form degree $5$. Consider a cohomology class $\alpha$ of the form
\[ \alpha_0 dx_1\cdots dx_5+ u \alpha_1 dx_1\cdots dx_5 + \cdots,\]
with $\alpha_i$'s in $\C[x_1,\ldots,x_5]\otimes_\C S$. The action of $\Z/5\Z$ fixes $dx_1\cdots dx_5$ and $u$. This implies that if $\alpha$ is a $\Z/5\Z$-invariant class, then the polynomials $\alpha_i$'s are also invariant. In other words, the polynomial degree (in the $x$'s variables) of $\alpha_i$'s are all divisible by $5$. This further implies for a homogeneous $\alpha$, we have
\[ \deg (\alpha) = \deg (\alpha_0) + \deg (dx_1\cdots dx_5) = -\frac{2}{5}(\mbox{ polynomial degree of $\alpha_0$})+3 \in 2\Z+1.\]
This shows that $HC^-(A^\mathcal{W})^{\Z/5\Z}$ is $\Z$-graded, in fact it is odd integer graded.
\end{proof}

We may split the group $(\Z/5\Z)^4\cong (\Z/5\Z)\oplus (\Z/5\Z)^3$ with $\Z/5\Z=\langle (1,1,1,1,1)\rangle$.  By Proposition~\ref{prop:equiv}, we have 
\[HC_{\sf odd}^-\big( A^\mathcal{W}\rtimes (\Z/5\Z)^4 \big)\cong \big( HC^-( A^\mathcal{W})\big)^{(\Z/5\Z)^4} \cong \Big(  \big( HC^-( A^\mathcal{W})\big)^{(\Z/5\Z)} \Big)^{(\Z/5\Z)^3}.\]
Since the VSHS $ \big( HC^-( A^\mathcal{W})\big)^{(\Z/5\Z)} $ is already $\Z$-graded, so is $HC_{\sf odd}^-\big( A^\mathcal{W}\rtimes (\Z/5\Z)^4 \big)$. Denote by $HC^-\big(A^\mathcal{W}\rtimes (\Z/5\Z)^4 \big)_{[1]}$ the equivalent description of VSHS as in Lemma~\ref{lem:equivalence}.

\begin{Theorem}~\label{thm:b-model-comparison}
With the notations set as above. There exists an isomorphism 
\[ \Phi:    HC^-\big(A^\mathcal{W}\rtimes (\Z/5\Z)^4 \big)_{[1]} \cong R^3\widetilde{\pi}_*\C\otimes_\C S,\]
as VSHS's over $B=\C\backslash\{1,e^{2\pi i/5}, e^{4\pi i/5}, e^{6\pi i/5}, e^{8\pi i/5}\}$.
\end{Theorem}

\begin{proof}
The idea of the proof is the following. The main work is to prove that there exists an isomorphism of VSHS's:
\[ HC^-\big( A^{\mathcal{W}}\rtimes (\Z/5\Z)\big)_{[1]} \cong R^3{\pi}_*\C\otimes_\C S.\]
Then we further take the $(\Z/5\Z)^3$-invariants on both side and use isomorphisms in~\ref{eq:iso1} and~\ref{eq:iso2} to deduce the result.

First, observe that there exists an isomorphism of VSHS's $$HC_{\sf odd}^-\big(A^{\mathcal{W}}\rtimes (\Z/5\Z)\big)\cong HC^-( A^\mathcal{W})^{\Z/5\Z},$$ proved in the same way as Proposition~\ref{prop:equiv}.  Taking the $\Z/5\Z$-invariants of $\Psi$ in Lemma~\ref{lem:comparison}, we obtain an isomorphism of $\Z$-graded VSHS's (still denoted by $\Psi$):
\[ \Psi: \big( HC^-(A^\mathcal{W})\big)^{\Z/5\Z} \cong H^*\big( \Omega_{B\times \C^5/B}^*[[u]], d\mathcal{W}+ud_{DR}\big)^{\Z/5\Z}.\]
To proceed from here, we construct a map
 \[\Theta: H^*\big( \Omega_{B\times \C^5/B}^*[[u]], d\mathcal{W}+ud_{DR}\big)^{\Z/5\Z}_{[1]}\ra  R^3\pi_*\C\otimes_\C S\]
via Griffiths' residue construction. Indeed, denote by ${{\sf Eu}}:= \sum_{j=1}^5 x_j \frac{\partial}{\partial x_j}$ the Euler vector field. Let $\alpha$ be a homogeneous element of the form $\alpha=[\alpha_0 dx_1\cdots dx_5+ u \alpha_1 dx_1\cdots dx_5 + \cdots]$, we define 
\begin{align*}
 \Theta(\alpha) &:= (-1)^{l_0}\cdot {{\sf Res}} \Big( \sum_i  (-1)^i (l_i-1)! \frac{\iota_{{\sf Eu}}( \alpha_i dx_1\cdots dx_5)}{\mathcal{W}^{l_i}}\Big)\\
 l_i &:= \frac{5-\deg(\alpha_idx_1\cdots dx_5)}{2}
 \end{align*}
 Note that by the previous discussion $l_i\in \Z$. The residue operator ${\sf Res}$ is defined in~\cite{Gri}. We need to verify that $\Theta$ vanishes on $(d\mathcal{W}+ud_{DR})$-exact terms. It suffices to show that for a homogeneous $4$-form $\beta\in \Omega^4_{B\times \C^5/B}$, we have that $\Theta(d\mathcal{W}\wedge \beta + ud_{DR})=0$. Indeed, set $l= \frac{5-\deg (d\mathcal{W}\wedge\beta)}{2}=\frac{6-\deg(\beta)}{2}$, we compute
 \begin{align*}
 &\Theta(d\mathcal{W}\wedge \beta + ud_{DR}\beta)\\
  = & {{\sf Res}} \Big((-1)^l (l-1)!\frac{\iota_{{\sf Eu}} (d\mathcal{W}\wedge \beta)}{\mathcal{W}^l}- (-1)^l(l-2)!\frac{\iota_{{\sf Eu}}(d_{DR}\beta)}{\mathcal{W}^{l-1}}\Big)\\
=&(-1)^l {{\sf Res}} \big( (l-1)!\frac{L_{{\sf Eu}} \mathcal{W}\wedge \beta -d\mathcal{W}\wedge \iota_{{\sf Eu}}\beta }{\mathcal{W}^l}  \big) - (-1)^l{{\sf Res}} \big( (l-2)! \frac{L_{{\sf Eu}}\beta-d_{DR}\iota_{{\sf Eu}}\beta}{\mathcal{W}^{l-1}}\big)
\end{align*}
We use $L_{\sf Eu}\mathcal{W}=5\mathcal{W}$, and furthermore, for the $4$-form $\beta$, its degree and its polynomial degree is related by $\deg(\beta) = - \frac{2}{5} ( \mbox{ polynomial degree of $\beta$} ) + 4$. Thus we obtain $L_{\sf Eu} \beta= \frac{20-5\deg(\beta)}{2}=5\cdot \frac{4-\deg(\beta)}{2}=5(l-1) $. Putting these two Lie derivatives in the above formula yields:
\begin{align*}
 &\Theta(d\mathcal{W}\wedge \beta + ud_{DR}\beta)\\
  = & (-1)^l {{\sf Res}} \Big( (l-2)!\frac{d_{DR}\iota_{{\sf Eu}}\beta}{\mathcal{W}^{l-1}}-(l-1)!\frac{d\mathcal{W}\wedge \iota_{{\sf Eu}}\beta }{\mathcal{W}^l}\big)\\
  = & (-1)^l {{\sf Res}} \Big( d_{DR}  [ (l-2)!\frac{\iota_{{\sf Eu}}\beta}{\mathcal{W}^{l-1}}] \Big)\\
  = & 0
\end{align*}
The last equality is because the residue map is defined in terms of integrals over cycles which vanish on exact forms. In~\cite{Gri}, Griffiths had already shown that $\Theta$ is an isomorphism, and that it respects the Hodge filtration. In the following, we verify that $\Theta$ also intertwines with the connection operators, and that it intertwines the pairing up to a constant factor. 

We first check the connection operator. Indeed, for $\alpha=[\alpha_0 dx_1\cdots dx_5+ u \alpha_1 dx_1\cdots dx_5 + \cdots]$, we have
\begin{align*}
& \Theta\nabla_{\frac{d}{d\psi}}(\alpha)\\
=& \Theta ( [\frac{d}{d\psi}\alpha]- [x_1x_2x_3x_4x_5\cdot \alpha] )\\
=& (-1)^{l_0}{{\sf Res}} \Big( \sum_i  (-1)^{i} (l_i-1)! \frac{\iota_{{\sf Eu}}( \frac{d}{d\psi}\alpha_i dx_1\cdots dx_5)}{\mathcal{W}^{l_i}}\Big) \\
&+(-1)^{l_0}{{\sf Res}} \Big( \sum_i  (-1)^i l_i! \frac{\iota_{{\sf Eu}}( x_1x_2x_3x_4x_5\alpha_i dx_1\cdots dx_5)}{\mathcal{W}^{l_i+1}}\Big)\\
= & (-1)^{l_0} {{\sf Res}} \Big( \sum_i  (-1)^{i} (l_i-1)! \frac{d}{d\psi} \big(\frac{\iota_{{\sf Eu}}( \alpha_i dx_1\cdots dx_5)}{\mathcal{W}^{l_i}}\big)\Big)\\
=& \nabla_{\frac{d}{d\psi}}\Theta(\alpha)
\end{align*}

Next, we check that $\Theta$ matches the pairing $(-,-)$ on $H^*\big( \Omega_{B\times \C^5/B}^*[[u]], d\mathcal{W}+ud_{DR}\big)^{\Z/5\Z}$ with the intersection pairing on $R^3\pi_*\C\otimes_\C S$, up to some univeral constant. This follows from Carlson-Griffiths' calculation of the intersection pairing through the Residue map~\cite[Theorem 3]{CarGri}. To state their result, let $fdx_1\ldots dx_5,gdx_1\ldots dx_5\in \big(\Omega^5_{B\times \C^5/B})^{\Z/5\Z}$ be two homogeneous $5$-forms, and denote by
\[a=-\frac{\deg(f)}{2}, \;\;\; b:=-\frac{\deg(g)}{2}.\]
With the degrees of $f$ and $g$ computed as by Equation~\ref{eq:degree}. Then Carlson-Griffiths' formula says that 
\[{\sf Res}( \frac{\iota_{\sf Eu} (fdx_1\ldots dx_5)}{\mathcal{W}^{a+1}})
\cap {\sf Res}( \frac{\iota_{\sf Eu} (gdx_1\ldots dx_5)}{\mathcal{W}^{b+1}})=
c_{a,b}\cdot {\sf Res}_0 \begin{bmatrix}
fgdx_1\ldots dx_5\\
\partial_1\mathcal{W},\ldots,\partial_5\mathcal{W}
\end{bmatrix} \]
where the constant $c_{a,b}:= (-1)^{a(a+1)/2+b(b+1)/2+3+b} \frac{5}{a!b!}$. From this formula, and that the pairing is only non-zero when $a+b=3$, we have
\begin{align*}
&\Theta(fdx_1\ldots dx_5) \cap \Theta(gdx_1\ldots dx_5) \\
=&  (-1)^{a+b}a!b!{\sf Res}( \frac{\iota_{\sf Eu} (fdx_1\ldots dx_5)}{\mathcal{W}^{a+1}})
\cap {\sf Res}( \frac{\iota_{\sf Eu} (gdx_1\ldots dx_5)}{\mathcal{W}^{b+1}})\\
=& (-1)^{a+b}a!b!c_{a,b}\cdot {\sf Res}_0 \begin{bmatrix}
fgdx_1\ldots dx_5\\
\partial_1\mathcal{W},\ldots,\partial_5\mathcal{W}
\end{bmatrix} \\
=& (-1)^{\frac{(a+b)(a+b+1)}{2}+(a+b)b}5\cdot {\sf Res}_0 \begin{bmatrix}
fgdx_1\ldots dx_5\\
\partial_1\mathcal{W},\ldots,\partial_5\mathcal{W}
\end{bmatrix} \\
=& (-1)^b 5\cdot {\sf Res}_0 \begin{bmatrix}
fgdx_1\ldots dx_5\\
\partial_1\mathcal{W},\ldots,\partial_5\mathcal{W}
\end{bmatrix} 
\end{align*}
On the other hand, by Lemma~\ref{lem:equivalence} and $\deg(fdx_1\ldots dx_5)=3-2a$, we have
\begin{align*}
&( fdx_1\ldots dx_5, gdx_1\ldots dx_5 ) \\
=& (\sqrt{-1})^{3-2a}\cdot {\sf Res}_0 \begin{bmatrix}
fgdx_1\ldots dx_5\\
\partial_1\mathcal{W},\ldots,\partial_5\mathcal{W}
\end{bmatrix} 
\end{align*}
Comparing the two identities yields that
\[ \Theta(fdx_1\ldots dx_5) \cap \Theta(gdx_1\ldots dx_5)= \sqrt{-1}\cdot 5\cdot (fdx_1\ldots dx_5, gdx_1\ldots dx_5 ).\]
In conclusion, by composition of $\Psi$ (after taking $\Z/5\Z$-invariants) and $\Theta$ we obtain an isomorphism
\[ \Theta\Psi: \big( HC^-(A^\mathcal{W})\big)^{\Z/5\Z}_{[1]}  \cong R^3\pi_*\C\otimes_\C S\]
which intertwines the connection operators, and preserves the pairing, up to a constant factor. From this, we deduce that the two VSHS's are isomorphic.  Further taking $(\Z/5\Z)^3$-invariants on both sides, we obtain an isomorphism of VSHS's:
\[  \big( HC^-(A^\mathcal{W})\big)_{[1]}^{(\Z/5\Z)^4} \cong \big( R^3\pi_*\C\otimes_\C S\big)^{(\Z/5\Z)^3} \]
Using the isomorphisms~\ref{eq:iso1}~\ref{eq:iso2}, the result follows.
\end{proof}

\begin{remark}
Ganatra-Perutz-Sheridan~\cite[Conjecture 1.14]{GPS} conjectured that for a smooth proper algebraic variety over a formal punctured disk, the VSHS associated with the derived category of coherent sheaves is isomorphic to Griffiths' VSHS. Theorem~\ref{thm:b-model-comparison} confirms this conjecture in the case of the mirror quintic family. To state this more precisely, consider the inclusion of rings
\[ i: S \hookrightarrow \mathbb{C}((q))=\mathbb{C}[q^{-1},q]], \;\; q:=\psi^{-1}.\]
Extension by scalar from $S$ to $\mathbb{C}((q))$ yields an isomorphism of VSHS's:
\[ \Phi:    HC^-\Big(\big(A^\mathcal{W}\otimes_S \mathbb{C}((q))\big)\rtimes (\Z/5\Z)^4 \Big)_{[1]} \cong R^3\widetilde{\pi}_* \mathbb{C}((q)).\]
Here we slightly abused the notation to still use  $\widetilde{\pi}: \widetilde{\mathfrak{X}}\ra {{\sf Spec\,}} \mathbb{C}((q))$ for the pull-back of $ \widetilde{\mathfrak{X}}\ra B={{\sf Spec\,}} S$ along the extension $i$.

The $A_\infty$-algebra $A^\mathcal{W}\otimes_S \mathbb{C}((q))$ may be seen as a minimal model of $\End ( \mathbb{C}((q))^{{\sf stab}})$, with $\mathbb{C}((q))^{{\sf stab}}$ a compact generator in the category of matrix factorizations $\MF(\mathcal{W})$ of $\mathcal{W}$ defined over the field $\mathbb{C}((q))$. Then, we  use Orlov's Calabi-Yau/Landau-Ginzburg correspondence~\cite{Orl} that there exists an equivalence
\[D^b({\sf coh}([\mathfrak{X}/(\Z/5\Z)^3]))\cong \MF_{(\Z/4\Z)^4}(\mathcal{W}),\]
both of which are $\mathbb{C}((q))$-linear categories. Thus the categorical VSHS's associated with these two categories are isomorphic. By Theorem~\ref{thm:b-model-comparison}, the right hand side VSHS indeed matches with the geometric VSHS, hence we deduce that the VSHS associated with the category $D^b({\sf coh}([\mathfrak{X}/(\Z/5\Z)^3]))$ is also isomorphic to $R^3\widetilde{\pi}_* \mathbb{C}((q))$, as conjectured by Ganatra-Perutz-Sheridan. 
\end{remark}

\paragraph{{\bf Calculation of the canonical primitive form.}} In this susbsection, we calculate the canonical primitive form associated with the category $\MF_{(Z/5Z)^4}(W)$, restricted to only the $\psi$-direction. More precisley, we shall again work with the formal setting, i.e. over the ring $\widehat{S}:=\mathbb{C}[[\psi]]$. In Theorem~\ref{thm:construction-vshs}, we obtained a primitive VSHS $\cV^{\MF_{(Z/5Z)^4}(W)}$, over a formal base $M={{\sf Spec}}\mathbb{C}[[HH^{[0]}\big(\MF_{(Z/5Z)^4}(W)\big)^\vee]]$. By Theorem~\ref{thm:canonical} we also obtain a canonical primitive form
\[ \zeta^{\MF_{(Z/5Z)^4}(W)} \in \cV^{\MF_{(Z/5Z)^4}(W)}.\]
The goal of this section is to compute $\zeta^{\MF_{(Z/5Z)^4}(W)}$, restricted to the marginal direction $\psi$, given by the dual coordinate of the Hochschild cohomology class $[x_1\cdots x_5]\in HH^{[0]}\big( \MF_{(Z/5Z)^4}(W)\big)$.

The maximal diagonal symmetry group of the Fermat quintic is $G_W^{{\sf max}}=(\Z/5\Z)^5$. We describe the canonical $G_W^{{\sf max}}$-equivariant splitting using the comparison isomorphism (Lemma~\ref{lem:comparison} and Proposition~\ref{prop:equiv}):
\[ HC^-_{\sf odd}\big(\MF_{(Z/5Z)^4}(W)\big) \cong H^*\big( \Omega_{\C^5}^*[[u]],dW+ud_{DR}\big)^{(\Z/5\Z)^4}.\]
In this model, the canonical splitting is simply given by the following basis vectors
\[ s_i:= (x_1x_2x_3x_4x_5)^i dx_1dx_2dx_3dx_4dx_5, \;\; i=0,1,2,3.\]
The flat extensions of $s_i$ along the $\psi$-direction can be computed by the following
\begin{Lemma}~\label{lem:trivial}
Let $\widehat{B}:={{\sf Spec \,}} \mathbb{C}[[\psi]]$. Then there exists an isomorphism of differential modules
\[  e^{(x_1x_2x_3x_4x_5)\psi/u}: H^*\big( \Omega_{\C^5}^*((u)),dW+ud_{DR}\big)\otimes_\C \widehat{S} \ra H^*\big( \Omega_{\C^5\times \widehat{B}/\widehat{B}}^*((u)),d\mathcal{W}+ud_{DR}\big).\]
Here the left hand side is endowed with the trivial connection $\frac{d}{d\psi}$, while the right hand side is endowed with the Gauss-Manin connection $\nabla^{{\sf GM}}_{\frac{d}{d\psi}}$.
\end{Lemma}

\medskip
Following a method of Li-Li-Saito~\cite{LLS}, we can compute the primitive form $\zeta$ associated with the splitting $s$ through its defining property that
\[ e^{\frac{(x_1x_2x_3x_4x_5)\psi}{u}}s_0= \zeta +\bigoplus_{k\geq 1} u^{-k} \widehat{S}\cdot e^{\frac{(x_1x_2x_3x_4x_5)\psi}{u}}\big({{\sf Im}} s\big).\]
Observe that we have
\begin{align*}
&(x_1x_2x_3x_4x_5)^n dx_1dx_2dx_3dx_4dx_5 \\
& = x_1^{n-4}(x_2x_3x_4x_5)^n dW\wedge dx_2\cdots dx_5\\
&= -u\cdot (n-4)\cdot x_1^{n-5}(x_2x_3x_4x_5)^n dx_1dx_2dx_3dx_4dx_5\\
&= u\cdot (n-4) \cdot x_1^{n-5} x_2^{n-4} (x_3x_4x_5)^n dW\wedge dx_1dx_3\cdots dx_5\\
&= u^2\cdot (n-4)^2\cdot (x_1x_2)^{n-5}(x_3x_4x_5)^n dx_1dx_2dx_3dx_4dx_5\\
&= \cdots\\
&= -u^5(n-4)^5(x_1x_2x_3x_4x_5)^{n-5} dx_1dx_2dx_3dx_4dx_5
\end{align*}

We take the differential $5$-form $$dx_1dx_2dx_3dx_4dx_5\in H^*\big( \Omega_{\widehat{B}\times \C^5/\widehat{B}}^*[[u]],d\mathcal{W}+ud_{DR}\big)^{(\Z/5\Z)^4},$$ pull it back through the isomorphism $ e^{\frac{(x_1x_2x_3x_4x_5)\psi}{u}}$ in Lemma~\ref{lem:trivial} and compute with the above formula to obtain
\begin{equation}~\label{eq:decomposition}
e^{-\frac{(x_1x_2x_3x_4x_5)\psi}{u}} dx_1dx_2dx_3dx_4dx_5= \omega_0(\psi)\cdot s_0-u^{-1}\omega_1(\psi)\cdot s_1+u^{-2}\omega_2(\psi)\cdot s_2 -u^{-3}\omega_3(\psi)\cdot s_3 
\end{equation}
with the power series $\omega_i(\psi)$ defined by
\[ \omega_i(\psi):= \sum_{k\geq 0} \frac{[(5k-4+i)(5k-9+i)\cdots (1+i)]^5}{(5k+i)!} \psi^{5k+i}.\]
Applying the flat extension operator $e^{\frac{(x_1x_2x_3x_4x_5)\psi}{u}}$  to Equation~\ref{eq:decomposition} above yields
\[ dx_1dx_2dx_3dx_4dx_5= \omega_0(\psi)\cdot s^{\sf flat}_0-u^{-1}\omega_1(\psi)\cdot s^{\sf flat}_1+u^{-2}\omega_2(\psi)\cdot s^{\sf flat}_2 -u^{-3}\omega_3(\psi)\cdot s^{\sf flat}_3\]
which shows that we have the following decomposition of $s_0^{\sf flat}$:
\[ s_0^{\sf flat}= \frac{1}{\omega_0}dx_1dx_2dx_3dx_4dx_5+ u^{-1}\frac{\omega_1}{\omega_0}\cdot s^{\sf flat}_1-u^{-2}\frac{\omega_2}{\omega_0}\cdot s^{\sf flat}_2 +u^{-3}\frac{\omega_3}{\omega_0}\cdot s^{\sf flat}_3\]
The primitive form $\zeta$ associated to a splitting, by the correspondence in Theorem~\ref{thm:bijection2}, is simply the postive part of $s^{\sf flat}_0$, and thus we deduce that
\[ \zeta= \frac{1}{\omega_0}\cdot dx_1dx_2dx_3dx_4dx_5.\]
This calculation proves the following LG/LG mirror symmetry result for the $B$-model orbifold $(W,(\Z/5\Z)^4)$ and its mirror dual $A$-model $(W,\Z/5\Z)$.
\begin{Theorem}~\label{thm:mirror}
Let $\cF_0$ denote the $g=0$ prepotential function associated with the category $\MF_{(\Z/5\Z)^4}(W)$ endowed with its canonical splitting. Then, the restriction of $\cF_0$ to the marginal deformation parameter in the flat coordinate $\tau:=\frac{\omega_1}{\omega_0}$ (or the mirror map) is equal to the genus zero FJRW prepotential function of $(W, \Z/5\Z)$ restricted to the marginal direction.
\end{Theorem}
\begin{proof}
Equation~\ref{eq:decomposition} matches with Chiodo-Ruan's calculation~\cite[Theorem 4.1.6]{ChiRua} of the $I$-function associated with the genus zero FJRW theory of $(W, \Z/5\Z)$. The procedure to obtain the prepotential function from the $I$-function is formal, see for example~\cite[Section 3]{ChiRua}.
\end{proof}

\end{document}